\magnification 1200
\font\BBfont=msbm10
\def\Bbb#1{\hbox {\BBfont #1}} 
\hsize15.5truecm
\baselineskip14pt
\font\myfont=cmr7 scaled \magstep1
%\font\sfont=cmbx7 scaled \magstep2
\font\lfont=cmbx10 scaled \magstep1
\font\bfont=cmbx10 scaled \magstep2
\def\litem{\par\noindent\hangindent=\parindent\ltextindent}

\def\ltextindent#1{\hbox to \hangindent{#1\hss}\ignorespaces}
\phantom{.}
\vskip 3truecm\noindent
{\centerline{\bfont {Ramsey and Nash-Williams combinatorics }}
\par\noindent
{\centerline{\bfont {via Schreier families}}}
 \vskip2pc\noindent
{\centerline {\bf Vassiliki Farmaki}}
\vskip2truecm\noindent
 {\centerline {\myfont{\bf (Preliminary Version)}}}
\vskip2truecm\noindent
 {\centerline {\myfont{\bf Abstract}}}
\vskip2pc
\par{\narrower\noindent
\par
{\baselineskip12pt {\myfont
The main results of this paper (a) extend the finite Ramsey partition theorem,  
and (b) employ this extension to obtain a stronger form of the infinite 
Nash-Williams partition theorem, and also a new proof   of Ellentuck's, 
and hence Galvin-Prikry's partition theorem. The proper tool for this 
unification of the classical partition theorems at a more general 
and stronger level is the system of Schreier families $({\cal A}_{\xi})$ 
of finite subsets of the set of natural numbers, defined for every countable ordinal $\xi$.
\par}}} 
\footnote{} {\myfont 2000 Mathematics Subject Classification: Primary 05D10;
Secondary 05C55.}
\footnote{} {\myfont Keywords and phrases: Ramsey, Nash-Williams combinatorics,  Schreier
families.}  
\vfill\eject
{\centerline{\lfont Introduction}}
\vskip1pc
The main results of this papers (Theorem A) extend the finite Ramsey ([R], 1929) partition theorem, and employ this extension to obtain a stronger form of the infinite Nash-Williams ([N-W], 1965) partition theorem (Theorems B, B$^\prime$, C), and also a new  proof of Ellentuck's ([E], 1974), and hence Galvin-Prikry's ([G-P], 1973) partition theorems.

Rather unexpectedly the ideas that lead to this purely combinatorial result have developed, in parallel, in Banach space theory, and they involve the system of Schreier families (${\cal A}_{\xi}$) of finite subsets of 
the set $\Bbb N$ of natural numbers defined for every countable ordinal $\xi$.
For the purpose of constructing a famous counterexample in the theory of Banach spaces, Schreier ([S], 1930) 
devised the classical Schreier family ${\cal F}_1$. This family definitely came into the attention of Banach space theory with the also famous Tsirelson  counterexample ([T], 1974), a construction that uses the Schreier family.  
In 1992 Alspach-Argyros ([A-A]) introduced the generalized Schreier families ${\cal F}_{\alpha}$ for every countable ordinal $\alpha$, and they used this family for the construction of Tsirelson type spaces, described by countable ordinals. In my 1994 paper ([F1]), working on some refinements of Rosenthal's ([R])paper on $c_0$ , 
it was  first realized that the generalized Schreier families ${\cal F}_{\alpha}$, corresponded to the ordinal 
$\omega^{\alpha}$, and that these was room for defining the intermediate generalized Schreier families 
${\cal A}_{\xi}$ for every countable ordinal $\xi$, thus filling the internals between 
$\omega^{\alpha}$ and $\omega^{\alpha+1}$, and in such a way that the earlier generalized Schreier families 
${\cal F}_{\alpha}$ essentially coincide with ${\cal A}_{\omega^{\alpha}}$. 
Independently, and working on a quite different problem (on distortable Banach spaces), Tomczak-Jaegermann 
([TJ], 1996), considered, at the suggestion of B. Maurey as she mentions, a variation of such families.

\par
At this point we might say that the complete system of thin Schreier families (${\cal A}_{\xi}$) (defined in 1.3 below), although a purely combinatorial object, arose in Banach space theory. There were some Banach space results that in retrospect can be considered as witnesses of a Ramsey type dichotomy nature of  these families or of their predecessors (Kiriakouli-Negrepontis ([M-N], 1992), Farmaki ([F1,F2], 1994), 
Argyros-Mercourakis-Tsarpalias ([A-M-T], 1998), Judd  ([J], 1999)).

In 1998 ([F4]) the study of the families (${\cal A}_{\xi}$) was refined  and employed for the proof of a far-reaching extension of the classical Ramsey theorem (Theorem A in this paper), one that holds for every countable ordinal $\xi$, of which the initial part, concerned with finite ordinals, coincides with the classical Ramsey theorem. 
Denoting by $[L]^{< \omega}$ the family of all finite subsets of a set $L$, here is the statement of the theorem:
\vfill\eject 
{\bf Theorem A (Ramsey partition theorem extended to countable ordinals)}. 
Let ${\cal F}$ be an arbitrary
family of finite subsets of $\Bbb N$, $M$ an infinite subset of $\Bbb N$ and
$\xi$ a countable ordinal number. Then, there exists an infinite subset $L$ of
$M$ such that 
$$ {\rm either} \ {\cal A}_{\xi} \cap [L]^{< \omega} \subseteq {\cal F} \ \ 
   {\rm or} \ \
 {\cal A}_{\xi} \cap [L]^{< \omega} \subseteq [\Bbb N]^{ < \omega}
  \setminus {\cal F} \ . $$
\vskip0.5pc\noindent
Since ${\cal A}_n = [\Bbb N]^{n}$ for any finite ordinal $n<\omega$, 
Theorem A (= Theorem 1.5 below) is in fact an extension of the classical Ramsey partition theorem from partitions on the families of n-tuples to (roughly speaking) partitions on the families ${\cal A}_{\xi}$ for any countable ordinal $\xi$.

The extended Ramsey Theorem A implies,  strengthened forms of the  Nash-Williams partition Theorem (Theorems B,B$^\prime$,C). We will employ Theorem A not for arbitrary, but only for hereditary families $\cal F$ of finite subsets of $\Bbb N$.  For such families ${\cal F}$ the strong Cantor-Bedixson index $s_L({\cal F})$ (Proposition 2.9) together with the canonical decomposition of any subset of 
$\Bbb N$ w.r.t. (${\cal A}_{\xi}$) (Proposition 2.4) imply a criterion that allows us to decide (in most cases) which horn of the dichotomy provided by Theorem A will actually hold.
Denoting by [L] the family of all infinite subsets of a set $L$, and by ${\cal A}_{\xi}^{\star}$ the family of all the initial segments of elements of ${\cal A}_{\xi}$, here is the statement of the theorem:
\vskip0.5pc
{\bf Theorem B 
(Stronger form of Nash-Williams partition theorem for hereditary families)}. 
Let $\cal F$ be a hereditary family of finite subsets of $\Bbb N$
and $M$ an infinite subset of $\Bbb N$. We have the
following cases:

{\bf [Case 1]} If the family ${\cal F}\cap [M]^{< \omega}$ is not pointwise closed,
then there exists $L \in [M]$ such that 
$[L]^{< \omega} \subseteq {\cal F}$.

{\bf  [Case 2]} If the family ${\cal F} \cap [M]^{< \omega}$ is pointwise  closed,
 then  setting
\par\noindent
{\centerline
{$ \xi^{\cal F}_M = sup \{s_L({\cal F}) : L \in [M] \} \ ;$}}
\par\noindent
which is a countable ordinal, the following subcases obtain:
\item{{\bf 2(i)} } For every countable ordinal $\xi$ with
 $\xi + 1 < \xi^{\cal F}_M$ there exists
$L \in [M]$ such that
\par\noindent
{\centerline{${\cal A}_{\xi} \cap [L]^{< \omega} \subseteq {\cal F} \ . $}}

\item{{\bf 2(ii)} } For every countable ordinal $\xi$
 with $\xi^{\cal F}_M < \xi +1$ there exists $L \in [M]$ such that
\par\noindent
{\centerline{$ {\cal F} \cap [L]^{< \omega} \subseteq ({\cal A}_{\xi})^{\star} 
   \setminus {\cal A}_{\xi} \ ;$}}
\item{} and equivalently,
\par\noindent
{\centerline{$ {\cal A}_{\xi} \cap [L]^{< \omega} \subseteq [\Bbb N]^{< \omega} 
    \setminus {\cal F} \ . $}}

\item{{\bf 2(iii)}} If $\xi^{\cal F}_M = \xi +1$, 
then both alternatives may materialize.
\vskip0.5pc
It is probably not apparent to the reader, why Theorem B (= Theorem 3.7 below) is in fact a result that deserves to be called a (strong) form of the Nash-Williams partition theorem. A convenient way to see this is by considering the reformulation that Gowers([G], 2002) gave of that theorem, and which can be stated as follows:
\vskip0.5pc
{\bf Nash-Williams partition theorem (in Gowers reformulation)}.
Let $\cal F$ be a family of finite subsets of $\Bbb N$. Then there exists an infinite subset $L$ of $\Bbb N$, such that 

either (i) $[L]^{<\omega} \subseteq {\cal F}$;
\par
or  \ \     (ii) for every infinite subsets $I$ of $L$, there exists an initial segment $s$ 
of $I$ which belongs to $[\Bbb N]^{< \omega} \setminus {\cal F}$.
\vskip0.5pc
Furthermore we remark that it is easy to see that 
WLOG we may assume in this reformulation that ${\cal F}$  be a {\bf tree} of finite subsets of $\Bbb N$
(cf. Remark 3.12).
 
A slightly weaker version of Theorem B, Theorem B$^\prime$, (= Theorem 3.10) concerns trees, and not necessarily hereditary families, of finite subsets of $\Bbb N$, bringing our result to a closer relation with the (tree form) Gowers reformulation of the Nash-Williams partition theorem.

Let us consider a further consequence of Theorem B$^\prime$ that brings forth in a clear manner the way in which our approach yields a result substantially stronger than the classical Nash-Williams result. The statement involves the decomposition (mentioned above) of any subset of $\Bbb N$ w.r.t. the system 
$({\cal A}_{\xi})_{\xi < \omega_1}$. In fact, by Proposition 2.4  below, every (infinite or finite) subset 
$I$ of $\Bbb N$ has a unique canonical representation w.r.t. each Schreier family ${\cal A}_{\xi}$, 
in such a way that for for every $\xi<\omega_1$, there is a unique initial segment  
  $s_{\xi,I}$ of $I$ that belongs to the Schreier family ${\cal A}_{\xi}$.
\vskip1pc
{\bf Theorem C (Stronger form of  Nash-Williams partition theorem in Gowers reformulation)}.
Let $\cal F$ be a tree of finite subsets of $\Bbb N$. Then there exists an infinite subset $L$ of 
$\Bbb N$, such that 

either (i) $[L]^{<\omega} \subseteq {\cal F}$;
\par
or  \ \     (ii)  there is a countable ordinal $\xi_0$, such that for every infinite subsets 
$I$ of $L$, there exists an initial segment $s$ of $I$ which belongs to $[\Bbb N]^{< \omega} \setminus {\cal F}$,
  and which is that unique initial segment of $I$ that belongs to ${\cal A}_{\xi_0}$.
\vskip0.5pc 
Compare this with the treeform  Gowers reformulation of the Nash-Williams theorem stated above. It is seen that our strengthened version provides, in the second horn of the dichotomy, not only the existence of the finite initial segments $s_I$ if $I$ (for all infinite subsets $I$ of $L$), but their determination by a countable ordinal 
$\xi_0$ in a {\bf unique} and {\bf uniform} way: thus the segment $s_I$ of $I$, does not simply exist as provided by the classical Nash-Williams result, but is that unique finite initial segment of $I$, which, according to the general decomposition of every subset of $\Bbb N$ w.r.t. the system $({\cal A}_{\xi})_{\xi< \omega_1}$, is an element of the family ${\cal A}_{\xi_0}$.

Ellentuck's theorem (not in a stronger form though), and thus the Galvin-Prikry partition theorem, follows also from our Theorem B (cf.  Theorem 4.6, Corollary 4.9, Remark 4.10).  

On the basis of these results, it is reasonable to conclude that the Schreier system 
$({\cal A}_{\xi})_{\xi < \omega_1}$, proves to be the correct combinatorial tool for the unification, extension and strengthening of all the finite Ramsey and infinite Nash-Williams partition theorems. These extended and strengthened partition theorems will no doubt find many applications, not only in the theory of Banach spaces, but in all the various areas where the classical combinatorial partition theorems have proved amply fruitful.
\vskip1pc
{\bf Notation.} We denote by $\Bbb N$ the set of all natural numbers. For an
infinite subset $M$ of $\Bbb N$ we denote by $[M]^{< \omega}$ the set of all
finite subsets of $M$, for $k \in \Bbb N$ we denote by $[M]^k$ the set of all
$k-$element subsets of $M$, and by $[M]$ the set of all infinite subsets of $M$ 
(considering them as strictly increasing sequences).

If $s,t$ are  non empty subsets of $\Bbb N$, then $s \preceq t$ means that $s$ is an
initial segment of $t$, while $s \prec t$ means that $s$ is a proper initial
segment of $t$. We write $s \leq t$ if $ max \ s \leq min \ t$,
while $s<t$ if $max \ s < min \ t$.

Identifying every subset of $\Bbb N$ with its characteristic function, we
topologize the set of all subsets of $\Bbb N$ by the topology of pointwise
convergence. 
\vskip3pc\noindent
{\lfont 1. The complete thin Schreier system and 
 the Ramsey 
 \par partition theorem extended to countable ordinals}
\vskip2pc
The main result in this section is a Ramsey type theorem for every countable ordinal $\xi$ 
(Theorem 1.5 - Theorem A), which can be considered as the countable ordinal analogue of the classical Ramsey theorem. This theorem is stated for the complete thin Schreier system ${(\cal A}_{\xi})_{\xi < \omega_1}$,
defined 1.3.

We recall Ramsey's classical partition theorem.
\vskip1pc
{\bf Theorem 1.1} ({\bf Ramsey} [R]). Let ${\cal F}$ be an arbitrary family of
finite subsets of $\Bbb N$, $M$ an infinite subset of $\Bbb N$ and $k$ a
natural number. Then there exists an infinite subset $L$ of $M$ such that 
$ {\rm either} \ \ [L]^k \subseteq {\cal F} \ \ \ {\rm or} \ \ \
[L]^k \subseteq [\Bbb N]^{< \omega} \setminus {\cal F} \ .$
\vskip0.5pc
This classical Ramsey partition theorem will prove to be, in Theorem A below,  the initial  segment 
of a whole family of
Ramsey type partition results,  one for every countable ordinal $\xi$.

In order to arrive at the statement of the  Ramsey partition theorem
for  any countable ordinal $\xi$ we need a $\xi$-ordinal
analogue ${\cal A}_{\xi}$ of $[\Bbb N]^k = {\cal A}_k$. This is accomplished for every 
$\xi < \omega_1$, by a rather laborious transfinite induction, that depends essentially
on a (classical) representation of (limit) ordinals, involving the ordinal
analogue of Euclidean algorithm as follows:
\vskip 1pc
{\bf Proposition 1.2} ({\bf Representation of ordinals}, [C2] [L]).
Let $\alpha$ be a non-zero, countable ordinal. For every limit ordinal $\xi$, so that
$ \omega^{\alpha} < \xi < \omega^{\alpha+1} \ , $
there exist a unique natural number $m \geq 0$, a  
sequence of ordinals $\alpha > \alpha_1 \ldots > \alpha_m > 0$ 
and natural numbers
$p,p_1, \ldots , p_m \geq 1$ (so that either $p>1$ or 
$p= 1$ and $m \geq 1$), such that
$ \xi = p\omega^{\alpha} + \sum^m_{i=1} p_i \omega^{\alpha_i} \ .$
\vskip1pc
We are now ready to define the families ${\cal A}_{\xi}$, for $\xi < \omega_1$,
which for reasons that will be explained later, will be collectively called the
complete thin Schreier system.
\vskip1pc
{\bf Definition 1.3} ({\bf The complete thin Schreier system 
$({\cal A}_{\xi})_{\xi<\omega_1 }$}).

For every non zero, limit ordinal $\lambda$ we choose and fix a strictly
increasing sequence $(\lambda_n)$ of successor ordinals smaller than $\lambda$
with ${\displaystyle sup_n} \ \lambda_n = \lambda$.
\par\noindent
We will define the system $({\cal A}_{\xi})_{\xi < \omega_1}$ recursively as follows:
\par\noindent
{\bf (1) [Case $\xi=0$]}\par
${\cal A}_0 = \{ \emptyset \} \ ; $
\par\noindent
{\bf (2) [Case $\xi = \zeta +1 $]}
\par\noindent
 $${\cal A}_{\xi} = {\cal A}_{\zeta +1} =
            \{ s \subseteq \Bbb N : s = \{ n \} \cup s_1 \ , \
{\rm where} \ 
          n \in \Bbb N , \{n \} < s_1 \
{\rm and} \ s_1 \in {\cal A}_{\zeta}\} \ ;$$
\par\noindent
{\bf (3) [Case $\xi = \omega^{\beta+1}, \ \beta$} countable ordinal{\bf]}
\par\noindent
$$\eqalign {
   {\cal A}_{\xi} = {\cal A}_{\omega^{\beta +1}} = {\cal B}_{\beta+1}=
   \{ s \subseteq \Bbb N : & \  s = \bigcup^n_{i=1} s_i \ {\rm with} \ 
                              n = min \ s_1 \ , \ \ s_1 < \ldots < s_n, \cr
                           & {\rm and} \ 
                          s_1, \ldots , s_n \in {\cal A}_{\omega^{\beta}} \};
  \cr} $$
\par\noindent
{\bf (4) [Case  $\xi =\omega^{\lambda}, \ \lambda$} non-zero, countable limit
ordinal {\bf]}\par
${\cal A}_{\xi} = {\cal A}_{\omega^{\lambda}} = {\cal B}_{\lambda}=
  \{ s \subseteq \Bbb N : s \in {\cal A}_{\omega^{\lambda_n}} \ {\rm with}
  \ \ n = min \ s \} \ , $
\par\noindent
(where $(\lambda_n)$ is the sequence of ordinals, converging 
to $\lambda \ $, fixed above); and
\par\noindent
{\bf (5) [Case $\xi \ {\rm limit} \ , \omega^{\alpha} < \xi < \omega^{\alpha +1}$} 
for some $0 < \alpha < \omega_1 ${\bf ]}

Let $\xi = p \omega^{\alpha} + {\displaystyle \sum^m_{i=1}} p_i
\omega^{\alpha_i} $ be the above representation (Proposition 1.2).
\par
$ {\cal A}_{\xi}= \{ s \subseteq \Bbb N : s = s_0 \cup 
   {\displaystyle (\bigcup^m_{i=1} s_i )} \ \ {\rm with} \ s_m < \ldots < s_1 < s_0 \ , $
\par $ \phantom {s_0 = }$
$ s_0 = s_1^0 \cup \ldots \cup s^0_{p} \ {\rm with} \ 
   s_1^0 < \ldots < s^0_p, \ s^0_j \in {\cal A}_{\omega^{\alpha}} , 
   1 \leq j
   \leq p \ , $
\par
$ s_i = s_1^i \cup \ldots \cup s_{p{_i}}^i, \ {\rm with} \
   s_1^i < \ldots < s_{p{_i}}^i , \ s_{j}^i \in 
   {\cal A}_{\omega^{\alpha_i}} , 1 \leq i \leq m , \ 1 \leq j \leq p_i \ \} \ . $
\vskip1pc
{\bf Remark 1.4} 
(i) ${\cal A}_{\xi} \subseteq [\Bbb N]^{< \omega} $ \ for every \ 
    $\xi < \omega_1 $ and $\emptyset \not\in {\cal A}_{\xi}$ for every $\xi > 0$.

(ii) ${\cal A}_k = [\Bbb N]^k $ \ \ for  \ \ $ k = 1,2, \ldots$

(iii) ${\cal B}_1 = {\cal A}_{\omega} = \{ s \in [\Bbb N]^{< \omega} :
       s = (n_1 < \ldots < n_k )$ \ \ with \ \ $n_1 = k \} \ . $
\par\noindent
Thus ${\cal A}_{\omega}$ is a modification of the {\bf classical 
 Schreier family}
([S]) 
\par
$
{\cal F}_1 = \{ s \subseteq \Bbb N : s = (n_1 < \ldots < n_k )
\ {\rm with} \ n_1 \geq k \} \ . $
\par\noindent
In this sense ${\cal A}_{\omega}$ is a thin Schreier family (this notion, used
also, by Pudlak - R\"odl, will be defined precisely later on in Definition 2.1).

(iv) ${\cal B}_{\alpha} = {\cal A}_{\omega^{\alpha}}$, for
$\alpha < \omega_1$, is defined using only
${\cal B}_{\beta} = {\cal A}_{\omega^{\beta}}$ for $\beta < \alpha $, and not
using all previously defined families 
${\cal A}_{\xi}, \ \xi < \omega^{\alpha}$. ${\cal B}_k$, for
$k \in \Bbb N$,   is a modification of 
generalized Schreier families defined by Alspach-Odell ([A-O]); and more
generally ${\cal B}_{\alpha}$, for $\alpha < \omega_1$ 
 is a modification of the families 
${\cal F}_{\alpha}$,  defined by Alspach-Argyros ([A-A]).
\vskip1pc
Now that the definition of the complete thin Schreier system
$({\cal A}_{\xi})_{\xi < \omega_1} $ is given, we are ready to state the first
 Ramsey partition theorem, for any countable ordinal $\xi$, a theorem whose
scope can be appreciated by the fact that the classical Ramsey theorem
corresponds to a finite ordinal $\xi < \omega$.
\vskip1pc
{\bf Theorem 1.5 (=Theorem A, Ramsey partition theorem extended to countable ordinals)}. 
Let ${\cal F}$ be an arbitrary
family of finite subsets of $\Bbb N$, $M$ an infinite subset of $\Bbb N$ and
$\xi$ a countable ordinal number. Then, there exists an infinite subset $L$ of
$M$ such that 
$$ {\rm either} \ {\cal A}_{\xi} \cap [L]^{< \omega} \subseteq {\cal F} \ \ 
   {\rm or} \ \
 {\cal A}_{\xi} \cap [L]^{< \omega} \subseteq [\Bbb N]^{ < \omega}
  \setminus {\cal F} \ . $$

In order to prove this theorem, we must find a way to relate the complete  thin
Schreier system  $({\cal A}_{\xi})_{\xi< \omega_1}$ with Ramsey type partition.
This is done in Proposition 1.7. The connecting concept that suits this purpose turns out to be the $\xi$-uniform families, which were defined by Pudl\'ak and R\"odl ([P-R]), an inductive concept incorporated in Proposition 1.7 below. 
\vskip1pc
{\bf Definition 1.6} Let  
${\cal L}$ be a family of finite subsets of $\Bbb N$. We set
\par
$ {\cal L}(n) = \{ s \in [\Bbb N]^{< \omega} : \{n \} < s \ \ {\rm and} \ \
   \{n \} \cup s \in {\cal L} \} $ for every $n \in \Bbb N$ .
 \vskip1pc  
{\bf Proposition 1.7}
For every countable ordinal $\xi$ there exists a concrete sequence $(\xi_n)$ of countable ordinals such that
$ {\cal A}_{\xi}(n) = {\cal A}_{\xi_n} \cap [\Bbb N \cap (n, + \infty)]^{< \omega} 
 {\rm for \ \ every} \ n \in \Bbb N $.
Moreover,  $\xi_n = \zeta$ for every $n \in \Bbb N$ if $\xi = \zeta +1$ is a successor ordinal and $(\xi_n)$ is a strictly increasing sequence with
$sup_n \xi_n = \xi$ if $\xi$ is a limit ordinal.
\vskip1pc
{\bf Proof}
We will prove it by recursion on $\xi$.
\par\noindent
{\bf (1) [Case $\xi =1$]} For every $n \in \Bbb N$ we have \par
$ {\cal A}_1(n)= \{ \emptyset \} = {\cal A}_0 \cap 
                   [\Bbb N \cap (n, + \infty)]^{< \omega} . $
\par\noindent
{\bf (2) [Case $\xi = \zeta +1$]} For every $n \in \Bbb N$ we have \par
$
{\cal A}_{\xi}(n)  = \{ s \subseteq \Bbb N : \{n\} < s 
                    \ \  {\rm and} \ s \in {\cal A}_{\zeta }\} = {\cal A}_{\zeta} 
                      \cap [\Bbb N \cap (n, + \infty)]^{< \omega} \ . $
\par\noindent
{\bf (3) [Case $\xi=\omega^{\beta+1}$ }for $0 \leq \beta < \omega_1  ${\bf ]}
For every $n \in \Bbb N$ we have 
$$ \eqalign {
{\cal A}_{\xi}(n) & = \{ s \subseteq \Bbb N : \{n \} < s, 
                      \{n \} \cup s = \bigcup^n_{i=1}  s_i \ ,
                      s_1 < \ldots < s_n \ {\rm and} \ s_1, \ldots , s_n \in 
                      {\cal A}_{\omega^{\beta}} \} \cr
            & =  \{ s \subseteq \Bbb N : s = s_0 \cup
                        ( \bigcup^n_{i=2} s_i) \ {\rm with} \
                        s_0 \in {\cal A}_{(\omega^{\beta})_n} \cap
                        [\Bbb N \cap (n, + \infty)]^{< \omega}, \cr
           & \ \ \ s_0<s_2< \ldots < s_n \ {\rm and} \ s_2, \ldots , s_n \in 
            {\cal A}_{\omega^{\beta}} \} \cr
                    & = {\cal A}_{ (n-1) \omega^{\beta} + 
                        (\omega^{\beta})_n } \cap 
                       [\Bbb N \cap (n, + \infty )]^{< \omega} \ , \cr} $$
according to the induction hypothesis.
Hence, \ \ $\xi_n = (n-1)\omega^{\beta} + (\omega^{\beta})_n$
\ for \  every \ \ $n \in \Bbb N$ \ and \ \ obviously 
${\displaystyle sup_n} \ \xi_n = \omega^{\beta+1} $. We note that in case 
$\xi = \omega$ we have $\xi_n = n-1$, since $\omega^0=1 \ . $ 
\par\noindent
{\bf (4) [Case $\xi = \omega^{\lambda}$ }for $\lambda$ non-zero, countable 
      limit ordinal{\bf ]} Let $(\lambda_n)$ be the sequence of successor
ordinals  converging to 
$\lambda $ fixed in the definition of the system
$({\cal A}_{\xi})_{{\xi}< \omega_1}$.
For every $n \in \Bbb N$ we have
$$\eqalign {
{\cal A}_{\xi}(n) & = \{ s \subseteq \Bbb N : \{n \}< s
                      \ \ {\rm and}\  \{n\} \cup s \in 
            {\cal A}_{\omega^{\lambda_n}} \} \cr 
         & = \{ s \subseteq \Bbb N : 
  s \in {\cal A}_{\omega^{\lambda_n}}(n) \} = 
 {\cal A}_{ (\omega^{\lambda_n})_n}  \cap 
  [\Bbb N \cap (n, + \infty)]^{< \omega} \ , \cr}  $$
according to the induction hypothesis.
If $\lambda_n = \mu_n +1$ for every
$n \in \Bbb N$, then  
$$ \xi_n =  (\omega^{\lambda_n})_n = 
 (\omega^{\mu_n+1})_n = (n-1) \omega^{\mu_n} + (\omega^{\mu_n})_n \ \
{\rm for \ every} \   n \in \Bbb N \ .$$
Of course ${\displaystyle sup_n} \ \xi_n = \omega^{\lambda}$,
since $\omega^{\mu_n} \leq (n-1) \omega^{\mu_n} \leq (n-1) 
\omega^{\mu_n} + (\omega^{\mu_n})_n \ .$
\par\noindent
{\bf (5) [Case $\xi$} limit, $\omega^{\alpha} < \xi < \omega^{\alpha +1}$
for some $0< \alpha < \omega_1$ {\bf ]} 
For simplicity we examine firstly the case  $\xi = p \omega^{\alpha}$ for some
$p \in \Bbb N$ with $p>1$.  
For every
$n \in \Bbb N$ we have
$$\eqalign { 
    {\cal A}_{\xi}(n) &  = \{ s \subseteq \Bbb N : \{ n \} < s \ {\rm and} \ \{n\} \cup s =
                           \bigcup_{i=1}^p s_i \ {\rm with} \ s_1 <\ldots< s_p \cr
                      & \ \ \ {\rm and} \ s_i \in {\cal A}_{\omega^{\alpha}} \ 
            {\rm for} \ 1 \leq i \leq p \}  \cr
                     & = \{s \subseteq \Bbb N : s = s_0 \cup 
                       (\bigcup^p_{i=2}s_i) \ {\rm with} \
                       \{n \} < s_0 < s_2< \ldots < s_p , \cr
                   & \ \ \ s_0 \in {\cal A}_{(\omega^{\alpha})_n} \ 
            {\rm and} \ s_2, \ldots , s_p 
                       \in {\cal A}_{\omega^{\alpha}} \} \cr
                   & = {\cal A}_{ (p-1) \omega^{\alpha} + (\omega^{\alpha})_n}
                       \cap [\Bbb N \cap (n, + \infty)]^{< \omega} \ . \cr} $$
\par\noindent
Thus \ \ $\xi_n = (p \omega^{\alpha})_n = (p-1)\omega^{\alpha} + (\omega^{\alpha})_n$
\ \ for \ every $n \in \Bbb N$. 
\ Of \  course, 
${\displaystyle sup_n} \ \xi_n= \xi$, since
${\displaystyle sup_n} \ (\omega^{\alpha})_n = \omega^{\alpha}$.

Now, let $\xi = \beta + p \omega^{\alpha}$,
 where $p \in \Bbb N$ with $p \geq 1$ and $\beta$ is an ordinal number with  
 $ 0< \beta < \xi$ (see Proposition 1.2). Then for every $n \in \Bbb N$ we have 
$$\eqalign {
{\cal A}_{\xi}(n) &  = \{s \subseteq \Bbb N :
                      \{n\} < s,  \{n\} \cup s = s_1 \cup s_2 ,
                   \ s_1 < s_2, \ s_1 \in {\cal A}_{p \omega^{\alpha}} \ 
                    {\rm and} \ s_2 \in {\cal A}_{\beta} \} \cr
                  & = \{ s \subseteq \Bbb N :
                      s = s_0 \cup s_2 \ {\rm with } \ \ \{n \} < s_0 <s_2 \ ,
                      s_0 \in {\cal A}_{p \omega^{\alpha}} (n) \ {\rm and} \ 
                     s_2 \in {\cal A}_{\beta} \} \cr
                  & = \{ s \subseteq \Bbb N : s = s_0 \cup s_2 \ 
                      {\rm with} \ \{n\} < s_0<s_2 ,
                     \ s_0 \in {\cal A}_{p \omega^{\alpha}}(n) \ 
                      {\rm and} \ s_2 \in {\cal A}_{\beta} \} \cr
                 & = {\cal A}_{\beta + (p \omega^{\alpha})_n }
                     \cap [\Bbb N \cap (n,+\infty)]^{< \omega} \ . \cr} $$
Hence, $\xi_n = \beta + (p \omega^{\alpha})_n =
        \beta + (p -1) \omega^{\alpha}+(\omega^{\alpha})_n $
for every $n \in \Bbb N$. Of course, ${\displaystyle sup_n} \ \xi_n = \xi$.
\par\noindent
This finishes the proof. 
\vskip0.5pc
We now, mimicking the standard proof of the classical Ramsey theorem, prove the $\xi$-Ramsey type theorem for every countable ordinal $\xi$.
\vskip1pc
{\bf Proof of Theorem 1.5}  
We will prove it by recursion on $\xi$.
Let $\xi =1$. Then,  ${\cal A}_1 = \{\{n\} : n \in \Bbb N \}$. Let ${\cal F} \subseteq [\Bbb N]^{< \omega}$ and 
$M \in [\Bbb N]$. Set $I = \{m \in M : \{m\} \in {\cal F} \}$ and consequently set $L=I$ in case $I$ is infinite
and $L= M \setminus I$ otherwise. Of course,
either ${\cal A}_1 \cap[L]^{< \omega} \subseteq {\cal F}$ or 
${\cal A}_1 \cap [L]^{< \omega} \subseteq [\Bbb N]^{< \omega} \setminus {\cal F}$.

Let $\xi > 1$. Assume that the theorem is valid for all ordinal $\zeta < \xi$. Let 
${\cal F} \subseteq [\Bbb N]^{<\omega}$ and $M \in [\Bbb N]$. Set $m_1 = min M$, $M_1 = M \setminus \{m_1\}$ and
${\cal F}_1 = \{ s \subseteq M_1 : \{ m_1\} \cup s \in {\cal F}\}$.  According to the previous proposition, 
${\cal A}_{\xi}(m_1) = {\cal A}_{\xi_{m_1}} \cap [\Bbb N \setminus \{1, \ldots , m_1 \}]^{< \omega}$
with $\xi_{m_1} < \xi$. Therefore, from the induction hypothesis, there exists $L_1 \in [M_1]$ such that
either ${\cal A}_{\xi_{m_1}} \cap [L_1]^{< \omega} \subseteq {\cal F}_1$ \ \ or \ \
${\cal A}_{\xi_{m_1}} \cap [L_1]^{< \omega} \subseteq [\Bbb N]^{< \omega} \setminus {\cal F}_1$.
Hence,
\par
either ${\cal A}_{\xi}(m_1) \cap [L_1]^{< \omega} \subseteq {\cal F}_1$ \ \ or \ \
${\cal A}_{\xi}(m_1) \cap [L_1]^{< \omega} \subseteq [\Bbb N]^{< \omega} \setminus {\cal F}_1$.
\par\noindent
Set $m_2 = min L_1$ ,  $M_2 = L_1 \setminus \{m_2\}$ and ${\cal F}_2 = \{s \subseteq M_2 : \{m_2\} \cup s \in {\cal F} \}$ .  Since
${\cal A}_{\xi}(m_2) = {\cal A}_{\xi_{m_2}} \cap [\Bbb N \setminus \{1, \ldots , m_2\}]^{< \omega}$ with 
$\xi_{m_2} < \xi$, and according to the induction hypothesis, there exists
$L_2 \in [M_2]$ such that

either ${\cal A}_{\xi} (m_2) \cap [L_2]^{< \omega} \subseteq {\cal F}_2$
or ${\cal A}_{\xi} (m_2) \cap [L_2]^{< \omega} \subseteq [\Bbb N]^{< \omega} \setminus {\cal F}_2$.
\par\noindent
Set $m_3 = min L_2$, $M_3 = L_2 \setminus \{m_3\}$ and proceed analogously.
\par\noindent
In this way we can construct a strictly increasing sequence 
$I=(m_n)_{n \in \Bbb N}$ in $M$, two decreasing sequences $(M_n)_{n\in \Bbb N}$,
$(L_n)_{n \in \Bbb N}$ in [M] such that:

(i) $m_{\kappa} \in L_n$ for every $\kappa > n$;

(ii) $L_n \subseteq M_n$ for every $n \in \Bbb N$; and

(iii) If ${\cal F}_n = \{ s \subseteq M_n : \{m_n\} \cup s \in {\cal F}\}$, then
\par\noindent
either ${\cal A}_{\xi} (m_n) \cap [L_n]^{< \omega} \subseteq {\cal F}_n$ \ \ or \ \ 
${\cal A}_{\xi} (m_n) \cap [L_n]^{< \omega} \subseteq [\Bbb N]^{< \omega} \setminus {\cal F}_n$.

Set $I_1 = \{m_n \in I $ : ${\cal A}_{\xi} (m_n) \cap [L_n]^{< \omega} \subseteq {\cal F}_n \}$
and $I_2 = I \setminus I_1$.
If $I_1$ is infinite, then ${\cal A}_{\xi} \cap [I_1]^{< \omega} \subseteq {\cal F}$.
Indeed, let $F \in {\cal A}_{\xi} \cap [I_1]^{< \omega}$ and let 
$m_n = min F$. Then $F = \{m_n\} \cup s$ for some $s \in {\cal A}_{\xi} (m_n)$.
 Since $m_n < s$ and $s \in [I]^{< \omega}$ we have that $s \subseteq L_n$ using (i).
Hence, $s \in {\cal A}_{\xi} (m_n) \cap [L_n]^{< \omega}$.
Since $m_n \in I_1$, we have that $s \in {\cal F}_n$, and consequently
$F = \{ m_n\} \cup s \in \cal F$.
\par\noindent
If $I_2$ is infinite, then, analogously can be proved that 
${\cal A}_{\xi} \cap [I_2]^{< \omega} \subseteq [\Bbb N]^{< \omega} \setminus \cal F$.
\par\noindent
Setting $L=I_1$ if $I_1$ is infinite, and $L = I \setminus I_1$ otherwise we have that

either
${\cal A}_{\xi} \cap [L]^{< \omega} \subseteq {\cal F}$ \ \ or \ \  
$ {\cal A}_{\xi} \cap [L]^{< \omega} \subseteq [\Bbb N]^{< \omega} \setminus {\cal F}$.

This finishes the proof.
\vskip1pc
{\bf Corollary 1.8} Let $M$ be an infinite subset of 
$\Bbb N, \{P_1, \ldots , P_n \}$ a finite partition of $[M]^{< \omega}$ and
$\xi$ a \ countable \ \ ordinal \ number. Then \ there \ exists $L \in [M]$ and 
$i \in \{1, \ldots, n \}$ such that
$ {\cal A}_{\xi} \cap [L]^{< \omega} \subseteq P_i $.
\vskip1pc
{\bf Corollary 1.9} Let ${\cal F}$ be a family of finite subsets of $\Bbb N$, 
$M \in [\Bbb N]$ and $\xi$ a countable ordinal. If ${\cal A}_{\xi} \cap 
{\cal F} \cap [I]^{< \omega} \ne \emptyset$ for every $I \in [M]$, then there
exists $L \in [M]$ such that
$ {\cal A}_{\xi} \cap [L]^{< \omega} \subseteq {\cal F} $.
\vskip3pc\noindent
{\lfont 2. Three basic properties of the complete thin Schreier system}
\vskip 2pc
In this section we will prove three basic properties of the complete thin Schreier systems 
$({\cal A}_{\xi})_{\xi < \omega_1}$, namely

(a) Each family ${\cal A}_{\xi}$ is thin (Proposition 2.2);

(b) every (finite of infinite) subset of $\Bbb N$ has (unique) canonical representation  
with respect to each family ${\cal A}_{\xi}$   (Proposition 2.4);
 and

(c) the strong Cantor-Bendixson index of ${\cal A}_{\xi}$ (given in
Definition 2.6 below) is precisely $\xi +1$; and this index is stable if
${\cal A}_{\xi}$ is restricted to 
${\cal A}_{\xi} \cap [M]^{< \omega}$ for any infinite subset $M$ of $\Bbb N$
(Proposition 2.9).

These properties will be necessary for establishing, Theorem  B in the next section 3, in case we have
not an arbitrary, but only a hereditary family $\cal F$ of finite subsets, an effective criterion that allows us to decide (is most cases) which horn of the dichotomy provided by Theorem A will actually hold. 
Although not immediately apparent, it turns out, as we shall see in section 3, that Theorem B is in fact a strengthened Nash-Williams partition theorem.

We start with some definitions.
\vskip1pc
{\bf Definition 2.1} Let ${\cal F}$ be a family of finite subsets of ${\Bbb N}$. 

(i) ${\cal F}$ is {\bf thin} if there are no elements $s,t \in {\cal F}$ 
with $s$ a proper initial segment (in the order of the natural numbers)  of $t$.

(ii) ${\cal F}^{\star} = \{t \in [\Bbb N]^{< \omega} : t$ is an initial
segment of some $s \in {\cal F} \}\cup \{ \emptyset \}$.

(iii)  ${\cal F}_{\star} = \{t \in [\Bbb N]^{< \omega} : t$ is a subset of some 
$s \in {\cal F} \} .$

(iv) ${\cal F}$ is {\bf hereditary} if ${\cal F}_{\star} = {\cal F}$.

(v) ${\cal F}$ is a {\bf tree} if ${\cal F}^{\star} = {\cal F}$.
\vskip1pc
{\bf Proposition 2.2}  Every family, ${\cal A}_{\xi}$ for $\xi < \omega_1$, is thin.
\vskip1pc
{\bf Proof}
We will prove it by induction on $\xi$.  The family ${\cal A}_1 = \{ \{n\} : n \in \Bbb N \}$
is obviously thin.
Let $\xi > 1$.  Assume that ${\cal A}_{\zeta}$ is thin, for every $\zeta < \xi$ .
Let $s \in {\cal A}_{\xi}$ and $t$ a proper initial semgent of $s$. If $t = \emptyset$, then
$t \not\in {\cal A}_{\xi}$ (Remark 1.4 (i)). Let $t \ne \emptyset$. If $m = mins$, then 
$t \setminus \{m\}$ is a proper initial segment of $s \setminus \{m\}$. Since 
$s \setminus \{m\} \in {\cal A}_{\xi}(m) \subseteq {\cal A}_{\xi_m}$ for some $\xi_m < \xi$ 
(Proposition 1.7), and ${\cal A}_{\xi_m}$ is thin, we have that 
$t \setminus \{m\} \not\in {\cal A}_{\xi_m}$ and consequently $t \not\in {\cal A}_{\xi}$.
This proves that ${\cal A}_{\xi}$ is a thin family.
\vskip1pc
In the following we will prove that every subset of $\Bbb N$ has 
 canonical representation with respect to each family ${\cal A}_{\xi}$. 
\vskip1pc
{\bf Definition 2.3} Let $\cal F$ be a family of finite subsets of $\Bbb N$. 

(i) A non-empty, finite subset $s$ of $\Bbb N$ has {\bf  canonical representation}
$R_{\cal F} (s) = \{ s_1, \ldots , s_n, s_{n+1} \}$ with type 
$t_{\cal F}(s)=n$ with respect to ${\cal F}$, if  there exist unique $n \in \Bbb N$,
$s_1, \ldots , s_n \in {\cal F}$ and $s_{n+1}$ a proper initial segment of some
element of $\cal F$ with $s_1< \ldots < s_n < s_{n+1}$ and such that
$s ={\displaystyle \bigcup^{n+1}_{i=1} } s_i \ .$

(ii) A infinite subset $I$ of $\Bbb N$ has {\bf  canonical representation}
$R_{\cal F} (I) = ( s_n)_{n \in {\Bbb N}}$
with respect to $\cal F$ if there exists unique sequence 
$(s_n)_{n \in \Bbb N}$ in $\cal F$ with 
$I = {\displaystyle \bigcup^{\infty}_{n=1}}s_n$ and such that $s_1 < s_2 < \ldots ,$
\vskip1pc
{\bf Proposition 2.4} ({\bf Canonical representation with respect
to ${\cal A}_{\xi}$ }) Let $\xi$ be a countable ordinal number. Every non-empty
subset of $\Bbb N$ has  canonical representation with respect to the
family ${\cal A}_{\xi}$.
\vskip1pc
{\bf Proof}
Let $I$ be an infinite subset of $\Bbb N$. We will prove, by induction on $\xi$, that $I$
 has canonical representation with respect to each family $\cal A_{\xi}$, $\xi < \omega_1$.
Of course $I$ has canonical representation with respect to the family ${\cal A}_1 = \{\{n\} : n \in \Bbb N\}$.
Let $\xi > 1$ and let $(\xi_n)$ be the corresponding sequence defined in Proposition 1.7.
Assume that the assertion holds for every $\zeta < \xi$.
Set $m_1 = min I$ and $I_1 = I \setminus \{m_1\}$. According to the induction hypothesis, $I_1$ has canonical representation $(s^1_n)_{n \in \Bbb N}$ with respect to 
${\cal A}_{\xi_{m_1}}$, since $1 \leq \xi_{m_1}< \xi$. Of course, 
$s_1 = \{m_1\} \cup s^1_1 \in {\cal A}_{\xi}$ .
Set $m_2 = min (I \setminus s_1)$ and $I_2 = I \setminus (s_1 \cup \{ m_2\})$.
According to the induction hypothesis, $I_2$ has canonical representation
$(s^2_n)_{n\in \Bbb N}$ with respect to
${\cal A}_{\xi_{m_2}}$,  since $1 \leq \xi_{m_2} < \xi$. Set $s_2 = \{ m_2\} \cup s^2_1$. 
Of course $s_2 \in {\cal A}_{\xi}$ and $s_1 < s_2$.
Set $m_3 = min (I \setminus ( s_1 \cup s_2))$ and $I_3 = I \setminus (s_1 \cup s_2 \cup \{ m_3\})$
and proceed analogously.
In this way, we can construct a sequence 
$(s_n)^{\infty}_{n=1}$ in ${\cal A}_{\xi}$ such that $s_1 < s_2 < \ldots$ and 
$I= {\displaystyle \bigcup^{\infty}_{n=1}} s_n$ This representation of $I$ with respect to $\cal A_{\xi}$
is unique, since $\cal A_{\xi}$ is a thin family (Proposition 2.2).
Hence, $I$ has canonical representation with respect to the family $\cal A_{\xi}$.

Now, let $\xi <\omega_1$ and $s = \{ m_1, \ldots , m_{\kappa}\}$,  be a non-empty, finite subset of 
$\Bbb N$ with $m_1 < \ldots < m_{\kappa}$. Set $m_{\kappa+i} = m_{\kappa}+i$ for every $i = 1,2, \ldots $ . 
The infinite set
$I =(m_n)^{\infty}_{n=1}$ has canonical representation $(s_n)^{\infty}_{n=1}$ with  respect to ${\cal A}_{\xi}$.
Using this fact, it is easy to prove that $s$ has canonical representation
with respect to $\cal A_{\xi}$.
\vskip1pc
{\bf Corollary 2.5} Let $\xi$ be a countable ordinal number. For every non-empty,
finite set $s$ of $\Bbb N$ exactly one of the following possibilities
occurs: 

either (i) $s$ is a proper initial segment of some element of ${\cal A}_{\xi}$;

or \ \ \ \ (ii) there exists an element of ${\cal A}_{\xi}$ which is an initial
segment of $s$.
\vskip1pc
{\bf Proof} Let $s \in [\Bbb N]^{< \omega}, s \ne \emptyset$. According to
Proposition 2.4, 
the case $t_{ {\cal A}_{\xi}}(s)=0$ gives
equivalently (i), while the complementary case, $t_{{\cal A}_{\xi}} (s) \geq 1$,
gives equivalently (ii).
\vskip1pc
In the following we will estimate the strong Cantor-Bendixson index of the families ${\cal A}_{\xi}$.
 This index (in Definition 2.6 below) is analogous to the
well-known Cantor-Bendixson index ([B],[C1]) and has been defined in [A-M-T].
Our notation is different from the one used by these authors.

We will prove in Proposition 2.9 below, that
the corresponding hereditary family of the thin Schreier family ${\cal A}_{\xi}$,
for $\xi < \omega_1$ has strong Cantor-Bendixson index equal to $\xi +1$, moreover for
every $M \in [\Bbb N]$ the restricted family
$({\cal A}_{\xi} \cap [M]^{< \omega})_{\star}$ has also index equal to $\xi +1 $.

This is the reason we have called $({\cal A}_{\xi})_{\xi < \omega_1}$ 
{\bf complete system.}
\vskip 1pc
{\bf Definition 2.6} ([A-M-T]) Let $\cal F$ be a hereditary and pointwise
closed family of finite subsets on $\Bbb N$. For $M \in [\Bbb N]$ we define the
{\bf strong Cantor-Bendixson  derivatives}
$({\cal F})^{\xi}_{M}$ of $\cal F$ on $M$ for every $\xi < \omega_1$ as follows:
$$ ({\cal F})^1_M = \{F \in {\cal F}[M]: F \ \
   {\rm is \ a \ cluster \ point  \ of } \ {\cal F}[F \cup L] \
   {\rm for \ each } \ L \in [M] \} ; $$
(where, ${\cal F}[M] = {\cal F} \cap [M]^{< \omega} )$.
\par\noindent
If $({\cal F})^{\xi}_M$ has been defined, then
$ ({\cal F})^{\xi +1}_M = \bigl( ({\cal F})^{\xi}_M \bigr )^1_M \ . $
\par\noindent
If $\xi$ is a limit ordinal and $({\cal F})^{\beta}_M$ have been defined for each
$\beta < \xi$, then

$ ({\cal F})^{\xi}_M = {\displaystyle \bigcap_{\beta < \xi}} 
   ({\cal F})^{\beta}_M \ . $

The {\bf strong Cantor-Bendixson index of ${\cal F}$ on $M$ } is defined to be the
smallest countable ordinal $\xi$ such that $({\cal F})^{\xi}_M = \emptyset$.
We denote this index by $s_M({\cal F})$.
\vskip1pc
{\bf Remark 2.7} (i) The strong Cantor-Bendixson index 
$s_M({\cal F})$ of a hereditary and
pointwise closed family $\cal F$ of finite subsets of $\Bbb N$ 
on some $M \in [\Bbb N]$
is a
countable successor  ordinal and  is less than  or equal to the ``usual"  
Cantor-Bendixson index
$O({\cal F})$ of $\cal F$ (see [K]).

\item{(ii)} If ${\cal F}_1,{\cal F}_2 \subseteq [\Bbb N]^{< \omega}$ are hereditary
and pointwise closed families
 with ${\cal F}_1 \subseteq {\cal F}_2$,  then 
$s_M({\cal F}_1) \leq s_M({\cal F}_2)$  for  every  
 $ M \in [\Bbb N]$ .

\item{(iii)} $s_M ({\cal F}) = s_M ({\cal F}\cap [M]^{< \omega})$ for every $M \in [\Bbb N]$.

\item{(iv)} For every $M \in [\Bbb N]$ and  $s \in [M]^{< \omega}$, 
according to a
remark in $[J]$, we have :
$$s \in ({\cal F})_M^1 \ {\rm if \ and \ only \ if \ the \ set} \ 
\{ m \in M : s \cup \{ m \} \not\in {\cal F} \} \ {\rm is \ finite} \ .$$

\item{(v)} \ Using the \ previous \ remark (iv), it can be \ proved by \ 
induction \ that \ for every 
$L \in [M]$ and  $\xi < \omega_1$
if $ \ A \in ({\cal F})^{\xi}_M, \ {\rm then} \ 
  F \cap L \in ({\cal F})^{\xi}_L $.
Hence, $s_L({\cal F}) \geq s_M ({\cal F})$. (see also [A-M-T]).

\item{(vi)} If $L$ is almost contained in $M$ (i.e. the relative complement 
$L \setminus M$ of $L$ in $M$ is a finite set), then 
$ s_L({\cal F}) \geq s_M({\cal F}) \ . $
\vskip1pc
In the following  we will calculate the strong Cantor-Bendixson index of the thin Schreier families 
${\cal A}_{\xi}$ .
\vskip1pc
{\bf Lemma 2.8} Let $\xi$ be a countable ordinal, $L \in [\Bbb N]$  and 
$\cal L$ a family of finite subsets of $\Bbb N$ such that ${\cal L}_{\star}$ 
and 
${\cal L}(n)_{\star} $ are pointwise closed for  every $n \in L$. 
\par
(i) If $F \in \bigl ( {\cal L}(n)_{\star} \bigr )^{\xi}_{L}$ for some
$n \in L$,  then 
$\{n\} \cup F \in ({\cal L}_{\star})^{\xi}_L$.
\par
(ii)  If
$F \ne \emptyset$ and $F \in ({\cal L}_{\star})^{\xi}_L$ , then there exist $ l \in \Bbb N$ with $l \leq min F$
and $I \in [L]$ such that $F \setminus \{l\} \in ({\cal L}(l)_{\star})^{\xi}_I$.
\vskip1pc
{\bf Proof} (i) We use induction on $\xi$. Let 
$F \in \bigl ({\cal L} (n)_{\star} \bigr )^1_L$.
Since
\par\noindent
$\{m \in M : F \cup \{m\} \in {\cal L}(n)_{\star} \} \subseteq
   \{m \in M : F \cup \{m\} \cup \{n\} \in {\cal L}_{\star} \} \ , $
we have, according to Remark 2.7 (iv), that 
$F \cup \{n\} \in ({\cal L}_{\star})^1_L \ . $
Let $1< \xi$. Suppose that the assertion holds for all ordinals $\zeta$ with $\zeta < \xi$.
If $F \in ({\cal L}(n)_{\star})^{\zeta +1}_L $,
then, according to the induction hypothesis,
$ \{ m \in M : F \cup \{ m \} \in ({\cal L}(n)_{\star})^{\zeta}_M \}
   \subseteq \{ m \in M : F \cup \{m\} \cup  \{n\} \in 
   ({\cal L}_{\star})^{\zeta}_M \} \ .$ 
Hence  
$\{n\} \cup F \in ({\cal L}_{\star})^{\zeta +1}_L$ (Remark 2.7 (iv)).
\par\noindent
The case where $\xi$ is a limit ordinal is trivial.

(ii) We use induction on $\xi$. Let $F \ne \emptyset$ and $F \in ({\cal
L}_{\star})^1_L $. According to Remark 2.7 (iv), the set 
$L_F = \{ m \in L : F \cup \{m\} \in {\cal L}_{\star}$ with $min F \leq m \}$
is almost equal to $L$. For each $m \in L_F$ there exists 
$s_m \in {\cal L}$ such that $F \cup \{m \} \subseteq s_m$. Of course
\par\noindent
$1 \leq min s_m \leq min F \ \ {\rm for \ every } \ m \in L_F \ $.  
 Set
\par
$    l = min \{n \in \Bbb N : \ {\rm the \ set} \ \{ m \in L_F :
   min \ s_m = n \} \ {\rm is \ infinite } \} ; \ {\rm and}  $
 \par
  $I = \{ m \in L_F : \ min \ s_m = l \} \cup F \ . $ 
\par\noindent
Then,  $\{l\} \leq F$,
$I \in [L] $ and $F \setminus \{l \} \in ({\cal L}(l)_{\star})^1_I$,
as required.

Suppose now that the assertion holds for all ordinals $\beta$ with 
$\beta < \xi$.  Firstly we examine the case 
$\xi = \zeta +1$. Let $F \ne \emptyset$ and $F \in ({\cal
L}_{\star})^{\zeta+1}_L$. According to Remark 2.7 (iv), the set
$L_F = \{ m \in L : F \cup \{m\} \in ({\cal L}_{\star})^{\zeta}_L$ and 
$min F \leq m \}$ is almost equal to $L$. Let $m_1 = min L_F$. By the induction
hypothesis there exist $l_1 \in \Bbb N$ with $l_1 \leq min F$ 
and $I_1 \in [L_F]$ such that 
$F \cup \{m_1 \} \setminus \{l_1\} \in 
 \bigl ({\cal L}(l_1)_{\star} \bigr )^{\zeta}_{I_1 \cup F}$, 
 since $F \cup \{ m_1 \} \in ({\cal L}_{\star})^{\zeta}_{L_F \cup F}$ (Remark 2.7 (v)). 
 Choose
$m_2 \in I_1$ and $m_2 > m_1$. Since 
$F \cup \{ m_2\} \in ({\cal L}_{\star})^{\zeta}_{I_1 \cup F}$,there 
\overfullrule=0ptexist
$l_2 \in \Bbb N$ with $l_2 \leq min F$ and $I_2 \in [I_1]$ such that
$F \cup \{m_2\} \setminus \{l_2\} \in \bigl ({\cal L}(l_2)_{\star}
 \bigr) ^{\zeta}_{I_2 \cup F} $. 
We continue analogously choosing $m_3 \in I_2$ with
$m_3 > m_2$ and so on.
Hence, we construct an increasing sequence 
$(m_i)^{\infty}_{i=1}$ in $L_F$, a sequence $(l_i)^{\infty}_{i=1}$ in $\Bbb N$,
with $1 \leq l_i \leq min F$ for every $i \in \Bbb N$, and a decreasing sequence
$(I_i)^{\infty}_{i=1}$ in $[L_F]$ such that 
$ F \cup \{ m_i \} \setminus \{l_i\} \in 
({\cal L}(l_i)_{\star})^{\zeta} _{I_i \cup F}$ for every $i \in \Bbb N$.
Let $l \in \Bbb N$ with $1 \leq l \leq  min F$ such that the set 
$L_1 =\{i \in \Bbb N : l_i = l \}$ is infinite. 
Set $I = \{ m_i : i \in L_1\} \cup F$. Then, 
$F \setminus \{l\} \in ({\cal L} (l)_{\star})^{\zeta+1}_I$, as required.

In the case where $\xi$ is a limit ordinal  we fix a strictly increasing sequence
$(\zeta_i)^{\infty}_{i=1}$ of ordinals with $\zeta_i < \xi$ for every 
$i \in \Bbb N$ and ${\displaystyle sup_i} \ \zeta_i = \xi$. 
Let $F \in ({\cal L}_{\star})^{\xi}_L$ and
$F \ne \emptyset$.
Then 
$F \in ({\cal L}_{\star})^{\zeta_i}_L$ for every $i \in \Bbb N$. 
According \ to
\ the \ induction \ \ hypothesis, 
there exist $l_1 \in \Bbb N$ with $l_1 \leq min F$ and  $I_1 \in [L \cap 
(min F, + \infty)]$
such that 
$F \setminus \{l_1\} \in ({\cal L}(l_1)_{\star})^{\zeta_1}_{I_1 \cup F}$.
Since $F \in ({\cal L}_{\star})^{\zeta_2}_{I_2 \cup F}$ there exists 
$l_2 \in \Bbb N$ with $l_2 \leq min F$ and $I_2 \in [I_1]$
such that $I_2 \ne I_1$ and
$ F \setminus \{l_2\} \in ({\cal L}(l_2)_{\star})^{\zeta_2}_{I_2 \cup F} \ .$
In this way, we construct a sequence $(l_i)^{\infty}_{i=1}$ with $1 \leq l_i \leq min F$
and a strictly decreasing sequence 
$(I_i)^{\infty}_{i=1}$ in $[L]$ such that
$
  F \setminus \{l_i\} \in ({\cal L}(l_i)_{\star})^{\zeta_i}_{I_i \cup F} \ , 
\ {\rm for \ every} \ i \in \Bbb N \ . $
\par\noindent
Let $l \in \Bbb N$  with $1 \leq l \leq min F$ such that the set
$L_1=\{i \in \Bbb N : l_i = l \}$ is infinite. Set
$I = \{ min I_i : i \in L_1 \} \cup F$. Then 
$F \setminus \{l\} \in ({\cal L}(l)_{\star})^{\zeta_i}_I$ for every
$i \in L_1$.Since ${\displaystyle sup_{i \in L_1} }\ \zeta_i = \xi $, we have that
$ F \setminus \{l\} \in \bigl (({\cal L}(l)_{\star} \bigr )^{\xi}_I \ . $

This completes the proof.
\vskip1pc
{\bf Proposition 2.9} {\bf (Cantor-Bendixson index of ${\cal A}_{\xi}$)}
Let $M$ be an infinite subset of $\Bbb N$. Then 
$$s_L(({\cal A}_{\xi} \cap [M]^{< \omega})_{\star}) = \xi +1 \ \ {\rm for \ every} \ \xi < \omega_1
\ \ \ {\rm and} \ L \in [M].$$

{\bf Proof}
It is easily proved, by induction on $\xi$, that the family $({\cal A}_{\xi} \cap [M]^{< \omega})_{\star}$
is pointwise  closed for every $\xi < \omega_1$ and $M \in [\Bbb N]$. Also, the family 
$(({\cal A}_{\xi} \cap [M]^{< \omega})(m))_{\star}$ is pointwise closed for every 
$\xi < \omega_1$, $M \in [\Bbb N]$ and $m \in M$, since, according to Proposition 1.7,
$({\cal A}_{\xi} \cap [M]^{< \omega})(m) =$ 
${\cal A}_{\xi_m} \cap [M \setminus \{ 1, \ldots , m\} ] ^{< \omega}$, where 
$\xi_m = \zeta$ if $\xi= \zeta+1$ is  a successor ordinal and $(\xi_m)_{m \in M}$ is a strictly increasing to $\xi$ sequence, if $\xi$ is a limit ordinal.

We will prove, by induction on $\xi$, that 
$s_L ( ({\cal A}_{\xi} \cap [M]^{< \omega})_{\star}) = \xi +1$ for every $\xi < \omega_1, M \in [\Bbb N]$ and 
$L \in [M]$.  Since $({\cal A}_1 \cap [M]^{< \omega})_{\star} = \{ \{ m\} : m \in M \} \cup \{ \emptyset\}$,
we have $(({\cal A}_1 \cap [M]^{< \omega})_{\star})^1_L = \{ \emptyset\}$ and consequently that
$s_L ( ({\cal A}_1 \cap [M]^{< \omega})_{\star}) = 2$ for every $M \in [\Bbb N]$ and $L \in [M]$.

Suppose that $\xi > 1$ and that the assertion holds for every ordinal $\zeta$ with $\zeta < \xi$. Let
$M \in [\Bbb N]$ and $L \in [M]$. For every $m \in M$, using Proposition 1.7, Remark 2.7(vi) and the induction hypothesis we get that $s_L ( ({\cal A}_{\xi} \cap [M]^{< \omega})(m)_{\star}) = \xi_m +1$ for every $m \in M$ and consequently that 
$\emptyset \in ( ({\cal A}_{\xi} \cap [M]^{< \omega})(m)_{\star})_L^{\xi_m}$ for every $m \in M$. In case $\xi = \zeta +1$ be a successor ordinal,
$\xi_m = \zeta$ for every $m \in M$ and, according to Lemma 2.8 (i),  
$\{l\} \in (({\cal A}_{\xi} \cap [M]^{< \omega})_{\star})^{\zeta}_L$ for every $l \in L$. 
Hence, $\emptyset \in (({\cal A}_{\xi} \cap [M]^{< \omega})_{\star})^{\xi}_L$ 
(Remark 2.7(iv)). On the other hand, in case $\xi$ be a limit ordinal, according to Lemma 2.8 (i),
$\emptyset \in (({\cal A}_{\xi} \cap [M]^{< \omega})_{\star})^{\xi_l}_L$ for every $l \in L$. Since 
$\xi_l < \xi$ and  ${\displaystyle sup_{l \in L} } \xi_l = \xi$,
we have also that $\emptyset \in (({\cal A}_{\xi} \cap [M]^{< \omega})_{\star})^{\xi}_L$.

In fact, $\{\emptyset\} = (({\cal A}_{\xi} \cap [M]^{< \omega})_{\star})^{\xi}_L$. Indeed, let 
$F \in (({\cal A}_{\xi} \cap [M]^{< \omega})_{\star})^{\xi}_L$ and $F \ne \emptyset$. Then, according to Lemma 2.8  (ii), there exist $n \in \Bbb N$ with $n \leq min F$ and $I \in [L]$ such that 
$F \setminus \{ n\} \in ( ({\cal A}_{\xi} \cap [M]^{< \omega})(n)_{\star})_I^{\xi}$.
This gives that 
$s_I ( ({\cal A}_{\xi} \cap [M]^{< \omega})(n)_{\star})\geq \xi +1 > \xi_n+1$.
A contradiction, according to Proposition 1.7 and the induction hypothesis.
Hence, $\{\emptyset\} = (({\cal A}_{\xi} \cap [M]^{< \omega})_{\star})^{\xi}_L$ and consequently 
$s_L(({\cal A}_{\xi} \cap [M]^{< \omega})_{\star})= \xi +1$ for every $\xi < \omega_1$.
\vskip 2pc\noindent
{\lfont 3. Strengthened Nash-Williams partition theorems}
\vskip2pc
We now turn our attention to the strengthened forms of the Nash-Williams theorem. These, contained in Theorems 3.7(= Theorem B), 3.10(= Theorem B$^{\prime}$), 3.11(= Theorem C) and 3.14, are consequences of the extended Ramsey Theorem 1.5 (= Theorem A), and of the tools contained in Section 2 (canonical representation 2.4, Cantor-Bendixson index 2.9)

Theorem 3.7(= Theorem B) can be considered as the extended Ramsey Theorem 1.5(= Theorem A), strengthened  for the case that we restrict ourselves, not to arbitrary, but only to hereditary families, ${\cal F}$ of finite subsets of $\Bbb N$.  As already remarked above it constitutes
in reality a strengthened Nash-Willimas type partition theorem, if we keep in mind the Gowers reformulation of Nash-Williams theorem (mentioned in the introduction above).
\vskip0.5pc
Proposition 3.1, a consequence of Theorem A, using also the canonical representation (Proposition 2.4), has consequences (Corollaries 3.4, 3.5) regarding generalized Schreier families (defined is 3.3).
\vskip1pc
{\bf Proposition 3.1} 
Let $\cal F$ be a family of finite subsets of 
$\Bbb N$ which is a tree 
\par\noindent
(${\cal F} = {\cal F}^{\star}$), $M$ an infinite subset of $\Bbb N$ and $\xi$ a countable ordinal number.
Then there exists $L \in [M]$ such that
\par\noindent
{\centerline{$
{\rm either} \ \ {\cal A}_{\xi} \cap [L]^{< \omega} \subseteq {\cal F} \ \ \
{\rm or} \ \ \ \ {\cal F} \cap [L]^{< \omega} \subseteq ({\cal A}_{\xi})^{\star}
 \setminus {\cal A}_{\xi} \ . $}}
\vskip0.5pc
{\bf Proof} According to the Ramsey partition  theorem  for the countable ordinal $\xi$
 there exists  $L \in [M]$ such that
$
{\rm either} \ \ {\cal A}_{\xi} \cap [L]^{< \omega} \subseteq {\cal F} \ \
{\rm or} \ \ {\cal A}_{\xi} \cap [L]^{< \omega} \subseteq [\Bbb N]^{< \omega}
 \setminus {\cal F} \ .$ 
Since $\cal F$ is a tree, we have
${\cal A}_{\xi} \cap [L]^{< \omega} \subseteq [\Bbb N]^{< \omega}
 \setminus {\cal F} $ if and only if 
 ${\cal F} \cap [L]^{< \omega} \subseteq ({\cal A}_{\xi})^{\star}
 \setminus {\cal A}_{\xi} \ . $
\par\noindent
Indeed, let ${\cal A}_{\xi} \cap [L]^{< \omega}
 \subseteq [\Bbb N]^{< \omega} \setminus {\cal F}$ and 
$F \in {\cal F} \cap [L]^{< \omega}$. According to Corollary 2.5, either there
exist $s \in {\cal A}_{\xi}$ such that $F$ is a proper initial segment of $s$
which gives that $F \in ({\cal A}_{\xi})^{\star} \setminus {\cal A}_{\xi}$,
as required, or there exists $t \in {\cal A}_{\xi}$ such that $t$ is an initial
segment of $F$. The second case is impossible. Indeed, since $\cal F$ is a tree
and $F \in {\cal F} \cap [L]^{< \omega}$, we have 
$t \in {\cal A}_{\xi} \cap [L]^{< \omega} \cap {\cal F}$. This 
contrary to our assumption that
${\cal A}_{\xi} \cap [L]^{< \omega} \subseteq [\Bbb N]^{< \omega} \setminus
 {\cal F}$. Hence, ${\cal F} \cap [L]^{< \omega} \subseteq ({\cal A}_{\xi})^{\star}
 \setminus {\cal A}_{\xi} $.
 \par\noindent
 It is obvious that if 
 ${\cal F} \cap [L]^{< \omega} \subseteq ({\cal A}_{\xi})^{\star}
 \setminus {\cal A}_{\xi} $, then 
 ${\cal A}_{\xi} \cap [L]^{< \omega} \subseteq [\Bbb N]^{< \omega}
 \setminus {\cal F} $.  
\vskip1pc
{\bf Corollary 3.2} Let $\xi_1, \xi_2$ be countable  ordinal numbers with
$\xi_1 < \xi_2$. For every $M \in [\Bbb N]$ there exists $L \in [M]$
such that
$({\cal A}_{\xi_1})_{\star} \cap [L]^{< \omega} \subseteq ({\cal A}_{\xi_2})^{\star}
  \setminus {\cal A}_{\xi_2} \ .$
\vskip1pc
{\bf Proof}  Of course 
$({\cal A}_{\xi_1})_{\star}$ is a tree. According to Proposition 3.1, 
 for every 
$M \in [\Bbb N]$ there exists $L \in [M]$ such that
\par
$
{\rm either} \ \ {\cal A}_{\xi_2} \cap [L]^{< \omega} \subseteq  ({\cal A}_{\xi_1})_{\star}\ \ \ 
{\rm or} \ \ \ ({\cal A}_{\xi_1})_{\star} \cap [L]^{< \omega} \subseteq ({\cal A}_{\xi_2})^{\star}
 \setminus {\cal A}_{\xi_2} \ . $
\par\noindent
The first alternative is impossible, since if  
${\cal A}_{\xi_2} \cap [L]^{< \omega}
 \subseteq ({\cal A}_{\xi_1})_{\star} \ $,
then
\par\noindent
$ \xi_2 +1 = s_L (({\cal A}_{\xi_2} \cap [L]^{< \omega})_{\star}) 
  \leq s_L (({\cal A}_{\xi_1})_{\star}) = \xi_1 +1 $ (Proposition 2.9).

A contradiction; hence 
$({\cal A}_{\xi_1})_{\star} \cap [L]^{< \omega}
  \subseteq ({\cal A}_{\xi_2})^{\star} \setminus {\cal A}_{\xi_2} \ . $
\vskip1pc
In the following, using Proposition 3.1, we indicate the close connection that
exists between the generalized Schreier families 
$({\cal F}_{\alpha})_{\alpha < \omega_1}$ and the $\omega^{\alpha}$- thin
Schreier families ${\cal A}_{\omega^{\alpha}}= {\cal B}_{\alpha}$ for 
$\alpha < \omega_1$. Firstly we will give
the appropriate definitions.
\vskip1pc
{\bf Definition 3.3} (i) (Generalized Schreier families
[S], [A-O],[A-A])

${\cal F}_0 = \{ \{n \} : n \in \Bbb N \} \cup \{ \emptyset \} \ ; $

${\cal F}_{\alpha +1} = \{F \subseteq \Bbb N : F = 
 {\displaystyle \bigcup^k_{i=1} } F_i \ , \{ k \} \leq
 F_1 < \ldots < F_k, \ {\rm and} \  F_i \in {\cal F}_{\alpha} \}
 \cup \{ \emptyset \} \ ; $
\par\noindent
If $\alpha$ is a limit ordinal choose and fix $(\alpha_n)_{n \in \Bbb N}$
strictly increasing to $\alpha$ and set

${\cal F}_{\alpha} = \{ F \subseteq \Bbb N : F \in {\cal F}_{\alpha_k} \ 
{\rm with} \  k \leq min F \} \cup \{ \emptyset \} \ . $

(ii) For a family ${\cal F}$ of finite subsets of $\Bbb N$ and 
$L = (l_n)^{\infty}_{n=1} \in [ \Bbb N]$ we set
\par\noindent
$ {\cal F}(L) = \{ (l_{n_1}, \ldots , l_{n_k})
   \in [L]^{< \omega} : (n_1, \ldots , n_k ) \in {\cal F} \} \ .$
\vskip1pc
{\bf Corollary 3.4} Let $\alpha$ be a countable ordinal. For every
$M \in [\Bbb N]$ there exists $L \in [M]$ such that
$ {\cal F}_{\alpha} (L) \subseteq ({\cal B}_{\alpha})^{\star}
   \subseteq {\cal F}_{\alpha} \ . $
\vskip1pc
{\bf Proof} The family ${\cal F}_{\alpha}$ is hereditary, hence, according to
Proposition 3.1, for every $M \in [\Bbb N]$ there exists $I \in [M]$ such that
\par\noindent
$
{\rm either} \ \ {\cal A}_{\omega^{\alpha}+1} \cap [I]^{< \omega}
  \subseteq {\cal F}_{\alpha} \ \ {\rm or} \ \ {\cal F}_{\alpha}
  \cap [I]^{< \omega} \subseteq ({\cal A}_{\omega^{\alpha}+1})^{\star}
  \setminus {\cal A}_{\omega^{\alpha}+1} \ . $
\par\noindent
The first alternative is impossible. 
Indeed, if 
${\cal A}_{\omega^{\alpha}+1} \cap [I]^{< \omega} \subseteq {\cal F}_{\alpha}$,
then \par\noindent
$ \omega^{\alpha} +2 = s_I ( ( {\cal A}_{\omega^{\alpha}+1}
   \cap [I]^{< \omega})_{\star}) \leq s_I ({\cal F}_{\alpha})=
   \omega^{\alpha}+1$ (Proposition 2.9).
\par\noindent
A contradiction; hence
${\cal F}_{\alpha} \cap [I]^{< \omega} \subseteq
 ({\cal A}_{\omega^{\alpha}+1})^{\star} \setminus
 {\cal A}_{\omega^{\alpha}+1}$.

Let $I = (i_n)^{\infty}_{n=1}$. We set 
$L=(i_n)^{\infty}_{n=3} = (l_n)^{\infty}_{n=1}$. We will prove that
\par\noindent
${\cal F}_{\alpha}(L) \subseteq ({\cal B}_{\alpha})^{\star}$. Indeed, let
$(l_{n_1}, \ldots , l_{n_k}) \in {\cal F}_{\alpha} (L) $, 
with  
$(n_1, \ldots , n_k) \in {\cal F}_{\alpha}$.  Then 
$(n_1+1, n_1+2, \ldots , n_k+2) \in {\cal F}_{\alpha}$ 
and consequently
$(i_{n_1 +1}, i_{n_1+2}, \ldots , i_{n_k +2}) \in
 {\cal F}_{\alpha} \cap [I]^{< \omega}$ (for the properties of 
${\cal F}_{\alpha}$ see [A-M-T]). This gives that
$(i_{n_1+1}, l_{n_1},\ldots l_{n_k}) \in 
 ({\cal A}_{\omega^{\alpha}+1})^{\star}$ and
consequently that $(l_{n_1}, \ldots , l_{n_k}) \in 
 ({\cal A}_{\omega^{\alpha}})^{\star} $ , as
required. Hence, ${\cal F}_{\alpha}(L) \subseteq ({\cal B}_{\alpha})^{\star}$.
It is obvious that $({\cal B}_{\alpha})^{\star} \subseteq {\cal F}_{\alpha}$.
  \vskip1pc
R. Judd in [J] had provided, using Schreier games, that for every
hereditary family $\cal F$ of finite subsets of $\Bbb N$, $\alpha <\omega_1$ and
$M\in [\Bbb N]$, either there exists $L\in [M]$ such that
${\cal F}_{\alpha}(L)\subseteq \cal F$ or there exists $L\in [M]$ and $N\in [\Bbb N]$
such that ${\cal F} \cap [N]^{< \omega}(L)\subseteq \cal F_{\alpha}.$

As a corollary of Proposition 3.1 we will prove a stronger version of 
this result.
\vskip1pc
{\bf Corollary 3.5} For every family $\cal F$ of finite subsets of $\Bbb N$ which is a tree, every
countable ordinal $\alpha$ and $M\in [\Bbb N]$ there exists $L\in [M]$ such that
\par\noindent
$
{\rm either} \ \ {\cal F}_{\alpha}(L)\subseteq {\cal F} \ \
 {\rm  or} \ \ {\cal F} \cap [L]^{<\omega}\subseteq
\cal F_{\alpha}\ . $
\vskip1pc
{\bf Proof} According to Proposition 3.1 there exists $N\in [M]$ such that
\par
$
{\rm either} \ \ {\cal B}_{\alpha} \cap [N]^{<\omega}\subseteq {\cal F} \ \ 
{\rm or} \ \ {\cal F} \cap [N]^{<\omega}\subseteq 
({\cal B}_{\alpha})^{\star} \ ; $
\par\noindent
If \ \  ${\cal B}_{\alpha} \cap [N]^{<\omega}\subseteq {\cal F}$, then, 
according to Corollary 3.4 and Proposition 2.4, there exists $L\in [N]$ 
such that
$
{\cal F}_{\alpha}(L)\subseteq ({\cal B}_{\alpha})^{\star}\cap [L]^{<\omega}
 \subseteq \bigl ({\cal B}_{\alpha} \cap [N]^
  {<\omega} \bigr)^{\star} \subseteq {\cal F} \ .$
\par\noindent
Hence, either ${\cal F}_{\alpha}(L)\subseteq {\cal F}$, 
or ${\cal F} \cap [L]^{<\omega}\subseteq ({\cal B}_{\alpha})^{\star}\subseteq 
  {\cal F}_{\alpha} \ .$ 
\vskip0.5pc
Since we will study the hereditary families of finite subsets of $\Bbb N$ which in addition are closed
in the pointwise topology we will
give an elementary characterization of them.
\vskip1pc
{\bf Proposition 3.6} Let $\cal F$ be a non empty,  family of finite
subsets of $\Bbb N$. 

(i) Let $\cal F$ be a tree. Then $\cal F$ is pointwise closed if and only if 
there does not exist an infinite sequence $(s_i)^{\infty}_{i=1}$ of
elements of $\cal F$ with $s_1 \prec s_2 \prec \ldots \ . $

(ii) Let $\cal F$ be hereditary. Then $\cal F$ is pointwise closed if and only if 
there does not exist $M \in [\Bbb N]$ such that
$[M]^{< \omega} \subseteq {\cal F}$.
\vskip1pc
{\bf Proof} (i)  Let $\cal F$ be a tree. If  $(s_i)^{\infty}_{i=1}$ is a sequence in
$\cal F$ with $s_1 \prec s_2 \prec \ldots$, then 
 $(s_i)^{\infty}_{i=1}$ converges
pointwise to an infinite subset $s$ of $\Bbb N$. Since
$s \not\in {\cal F}$, $\cal F$ is not closed.
 We assume that  there does not exist an infinite sequence $(s_i)^{\infty}_{i=1}$ of
elements of $\cal F$ with $s_1 \prec s_2 \prec \ldots $.  Let
 $(t_n)^{\infty}_{n=1} \subseteq {\cal F}$
converges pointwise to some subset $t$ of $\Bbb N$.
If $t$ is finite, then $t$ is an initial segment of some $t_{n_0}$  for some $n_0 \in \Bbb N$. 
Since $\cal F$ is a tree, $t \in \cal F$ .
If $t = (n_1, n_2, \ldots)$ with $n_1<n_2 < \ldots$,  then we
set $s_i =(n_1,n_2, \ldots, n_i)$ for every $i \in \Bbb N$.
Of course  $s_1 \prec s_2 \prec \ldots$. Let $s_n^i = t_n \cap [0,n_i]$ 
 for every $i \in \Bbb N$ and $n \in \Bbb N$.
It is easy to see that the sequence 
$(s_n^i)^{\infty}_{n=1}$ in $\cal F$ converges pointwise to $s_i$. 
 According to the previous  case,   $s_i \in \cal F$,
for every $i \in \Bbb N$. A contradiction to our assumption, 
 so $t$ is finite and  $t \in {\cal F}$.  Hence, $\cal F$ is pointwise closed.

(ii) It is easily proved, using (i).
\vskip1pc
Now, using Propositions 3.1, 3.6 and  the concept of the strong Cantor-Bendixson index 
(Proposition 2.9) we state and prove the stronger form of the Nash-Williams partition theorem for hereditary families of finite subsets of $\Bbb N$.
\vskip1pc
{\bf Theorem 3.7 (=Theorem B, Stronger form of Nash-Williams partition theorem for hereditary families)} 
Let $\cal F$ be a hereditary family of finite subsets of $\Bbb N$
and $M$ an infinite subset of $\Bbb N$. We have the
following cases:
\vskip0.5pc\noindent
{\bf [Case 1]} The family ${\cal F}\cap [M]^{< \omega}$ is not pointwise closed.
\par\noindent
Then, there exists $L \in [M]$ such that 
$[L]^{< \omega} \subseteq {\cal F}$.  
\vskip0.5pc\noindent
{\bf  [Case 2]} The family ${\cal F} \cap [M]^{< \omega}$ is pointwise  closed.
Then, setting
$$ \xi^{\cal F}_M = sup \{s_L({\cal F}) : L \in [M] \} \ ,$$
which is a countable ordinal, the following subcases obtain:
\par\noindent
{\bf 2(i)} If $\xi^{\cal F}_M > \xi +1  $, then there exists $L \in [M]$ such that
$$ {\cal A}_{\xi} \cap [L]^{< \omega} \subseteq {\cal F} \ ; $$
\par\noindent
{\bf 2(ii)} if $\xi^{\cal F}_M < \xi +1$ , then for every $I \in [M]$ there exists $L \in [I]$ 
$$ {\cal A}_{\xi} \cap [L]^{< \omega} \subseteq [\Bbb N]^{< \omega} 
   \setminus {\cal F} \ ; \ ({\rm \ equivalently \ such \ that}
 {\cal F} \cap [L]^{< \omega} \subseteq ({\cal A}_{\xi})^{\star}
   \setminus {\cal A}_{\xi}) \ ;\ {\rm and}, $$
\par\noindent
{\bf 2(iii)} if $\xi_M^{\cal F} = \xi+1$, then there exists $L \in [M]$ such that
$$
{\rm either} \ \ {\cal A}_{\xi} \cap [L]^{< \omega} \subseteq {\cal F} \ \
{\rm or} \ \ {\cal A}_{\xi} \cap [L]^{< \omega} \subseteq [\Bbb N]^{< \omega} \setminus {\cal F} .$$
\par\noindent
Both alternatives in 2(iii) may materialize.
\vskip1pc
{\bf Proof } [Case 1] If the family ${\cal F} \cap [M]^{< \omega}$ is not pointwise closed,
then there exists $L \in [M]$ such that 
$[L]^{< \omega} \subseteq {\cal F}$,  according to Proposition 3.6.
\par\noindent
[Case 2] Let ${\cal F} \cap [M]^{< \omega}$ be pointwise closed. Then $\xi^{\cal F}_{M}$ 
is a  countable ordinal. Indeed,
since the Cantor-Bendixson index $O({\cal F})$
of $\cal F$ (see [K]) is a countable ordinal
(as the family of
derived sets of $\cal F$ is countable) and since
$ s_I ({\cal F})\leq O({\cal F}) \ \ {\rm for \ every} \ \ I \in [\Bbb N] \ , $
we have  $\xi^{\cal F}_M \leq O ({\cal F}) < \omega_1$. 

2(i) Let $\xi +1 < \xi^{\cal F}_M$. Then, there exists $I \in [M]$ such that
$\xi+1 < s_I ({\cal F})$.
According to Theorem 1.5 and 
Proposition 3.1, there exists $L \in [I]$ such that
\par\noindent
$
{\rm either} \ \ {\cal A}_{\xi} \cap [L]^{< \omega} \subseteq {\cal F} \ \
{\rm  or} \ \  {\cal F} \cap [L]^{< \omega} \subseteq ({\cal A}_{\xi})^{\star}
   \setminus {\cal A}_{\xi} \subseteq ({\cal A}_{\xi})_{\star} \ . $
\par\noindent
The second alternative is impossible. Indeed, if
  ${\cal F} \cap [L]^{< \omega} \subseteq ({\cal A}_{\xi})_{\star}$, 
then, using Proposition 2.9 and Remark 2.7, we have
\par\noindent
$ s_I({\cal F}) \leq s_L ({\cal F}) = s_L ({\cal F} \cap [L]^{< \omega})
   \leq s_L ( ({\cal A}_{\xi})_{\star}) = \xi +1 \ . $
\par\noindent
This is a  contradiction;
hence ${\cal A}_{\xi} \cap [L]^{< \omega} \subseteq {\cal F} \ . $

2(ii) Let $\xi^{\cal F}_M < \xi +1$ and $I \in [M]$. 
 According to the Ramsey partition type theorem for the countable ordinal $\xi$
(Theorem 1.5), there exists 
$L \in [I]$ such that
\par\noindent
$
{\rm either}\ \ {\cal A}_{\xi} \cap [L]^{< \omega} \subseteq {\cal F} \ \
{\rm  or} \ \ 
       {\cal A}_{\xi} \cap [L]^{< \omega} \subseteq [\Bbb N]^{< \omega}
         \setminus {\cal F} \ . $
\par\noindent
The first alternative is impossible. Indeed, 
if ${\cal A}_{\xi} \cap [L]^{< \omega} \subseteq {\cal F}$,
then,  using  Proposition 2.9 and Remark 2.7,
we obtain 
$ \xi +1  = s_L (({\cal A}_{\xi} \cap [L]^{< \omega})_{\star})
   \leq s_L ({\cal F}) \leq \xi^{\cal F}_M \ . $
\par\noindent
This is a contradiction; hence, 
$ {\cal A}_{\xi} \cap [L]^{< \omega} \subseteq [\Bbb N]^{< \omega} 
  \setminus {\cal F} $ and
equivalently, 
\par\noindent
${\cal F} \cap [L]^{< \omega} \subseteq ({\cal A}_{\xi})^{\star}
 \setminus {\cal A}_{\xi} $, according to Proposition 3.1.

2(iii) That both alternatives in the case $\xi^{\cal F}_{M} = \xi +1$ may materialize can be seen 
by considering  two simple examples:

(1)
$ {\cal F} = \{ s \in [\Bbb N]^{< \omega} : s \ne \emptyset \ {\rm and} \
   \vert s \vert = 2 min \ s +1 \} , $
\par\noindent
where $\vert s \vert$ denotes the cardinality of $s$.
\par\noindent
(It is easy to see that 
$ {\cal F}(n) = [ \Bbb N \cap (n, + \infty)]^{2n} = {\cal A}_{2n} \cap [ \Bbb N \cap (n, + \infty)]^{< \omega}$
for every $n \in \Bbb N$.
The family  ${\cal F}_{\star}$ is pointwise closed and according to Lemma 2.8,
$s_I({\cal F}_{\star}) = \omega +1$ for every $I \in [\Bbb N]$.
Hence 
$ \xi^{{\cal F}_{\star}}_{M}= \omega+1 $ for every $M \in [\Bbb N]$.
It is now easy to verify that
\par\noindent
${\cal A}_{\omega} \cap [L]^{< \omega} \subseteq {\cal F}_{\star} \  
 {\rm for \ every} \  L \in [M]) \ ;$ and, 

(2) $ {\cal F} = \{ s \in [M]^{< \omega} : s \ne \emptyset \ {\rm and} \
   \vert s \vert = {{ min \ s} \over 2} \} \ , $
\par\noindent
  where $M$ stands for all  non zero, even natural numbers.
\par\noindent
(Since ${\cal F}(m) = {\cal A}_{ {m \over 2}-1} \cap [ M \cap (m,+ \infty)]^{< \omega}$
for every $m \in M$,
from Lemma 2.8 we get that $s_I ({\cal F}_{\star}) = \omega +1$ for every 
$I \in [M]$.
Thus 
$\xi^{{\cal F}_{\star}}_M = \omega+1 \ .$
It is now easy to verify that
$ {\cal F}_{\star} \cap [L]^{< \omega} \subseteq ({\cal A}_{\omega})^{\star} \setminus
   {\cal A}_{\omega} \ \ {\rm for \ every} \ \ L \in [M] $ and, according to Proposition 2.1, 
that ${\cal A}_{\omega}\cap [L]^{< \omega} \subseteq [\Bbb N]^{< \omega} \setminus {\cal F}_{\star}$ 
for every $L \in [M]$.)
\vskip1pc
As a corollary of Theorem 3.7 we have the following result of Argyros,
Merkourakis and Tsarpalias ([A-M-T]). 
\vskip1pc
{\bf Corollary 3.8} Let ${\cal F}$ be a hereditary and pointwise closed family
of finite subsets of $\Bbb N$, If there exists $M \in [\Bbb N]$ such that 
$s_M ({\cal F}) \geq \omega^{\alpha}$, then there exists $L \in [M]$
such that ${\cal F}_{\alpha} (L) \subseteq {\cal F}$.
\vskip1pc
{\bf Proof} If $s_M({\cal F}) > \omega^{\alpha} +1$, then, according to Theorem
3.7, there exists $N \in [M]$ such that
 ${\cal B}_{\alpha} \cap [N]^{<\omega}\subseteq \cal F$ and
according to Corollary 3.4 and Proposition 2.4 there exists $L\in [N]$ such that
${\cal F}_{\alpha}(L)\subseteq ({\cal B}_{\alpha})^{\star}\cap [L]^{<\omega}
\subseteq ({\cal B}_{\alpha} \cap [N]^{< \omega})^{\star}\subseteq {\cal F} $.

Now, if $s_M({\cal F})=\omega^{\alpha}+1$, then we set 
${\cal L}=\{\{m\}\cup s:s\in {\cal F}, m \in M \ {\rm and} \  \{m\}<s\} \ .$
It is easy to see that
$s_M({\cal L})>\omega^{\alpha}+1$. 
So applying the previous case to
the family $\cal L$ we can find $N= (n_i)^{\infty}_{i=1}\in [M]$ such that
${\cal F}_{\alpha}(N)\subseteq {\cal L}$. Setting
$L=(n_i)^{\infty}_{i=3}$ we have that ${\cal F}_{\alpha}(L)\subseteq \cal F$,
 as required. 
 \vskip1pc
 The version of Theorem B for trees is given below.
\vskip1pc
{\bf Definition 3.9} Let $\cal F$ be a family of finite subsets of $\Bbb N$. We set 

(i) ${\cal F}_h = \{ s \in \cal F$ : every non- empty subset of $s$ belongs to 
$\cal F \} \cup \{\emptyset \}$, and
\par\noindent
Of course, ${\cal F}_h$ is the largest subfamily of $\cal F$ which is hereditary.
\vskip1pc
{\bf Theorem 3.10  (= Theorem B$^{\prime}$, Stronger form of Nash-Williams partition theorem for trees)} 
Let $\cal F$ be a tree  of finite subsets of $\Bbb N$
and $M$ an infinite subset of $\Bbb N$. We have the
following cases: 
\vskip0.5pc\noindent
{\bf [Case 1]} The family ${\cal F}_h\cap [M]^{< \omega}$ is not pointwise closed.
\par\noindent
Then, there exists $L \in [M]$ such that 
$[L]^{< \omega} \subseteq {\cal F}$.  
\vskip0.5pc\noindent
{\bf  [Case 2]} The family ${\cal F}_h \cap [M]^{< \omega}$ is pointwise  closed.
Then setting
$$ \zeta^{\cal F}_M = sup \{s_L({\cal F}_h) : L \in [M] \} = \xi_M^{{\cal F}_h}\ ,$$
which is a countable ordinal, the following subcases obtain:
\par\noindent
{\bf 2(i)} If $\zeta^{\cal F}_M > \xi +1  $, then there exists $L \in [M]$ such that
$$ {\cal A}_{\xi} \cap [L]^{< \omega} \subseteq {\cal F} \ ; $$
\par\noindent
{\bf 2(ii)} if $\zeta^{\cal F}_M < \xi $ , then for every $I \in [M]$ there exists $L \in [I]$ such that
$$ {\cal A}_{\xi} \cap [L]^{< \omega} \subseteq [\Bbb N]^{< \omega} 
   \setminus {\cal F} \ ; \ ({\rm \ equivalently} \ 
 {\cal F} \cap [L]^{< \omega} \subseteq ({\cal A}_{\xi})^{\star}
   \setminus {\cal A}_{\xi}) \ ;\ {\rm and}, $$
\par\noindent
{\bf 2(iii)} if $\zeta_M^{\cal F} = \xi+1$ or $\zeta_M^{\cal F} = \xi$, then there exists $L \in [M]$ such that
$$
{\rm either} \ \ {\cal A}_{\xi} \cap [L]^{< \omega} \subseteq {\cal F} \ \
{\rm or} \ \ {\cal A}_{\xi} \cap [L]^{< \omega} \subseteq  [\Bbb N]^{< \omega} \setminus {\cal F} \ . $$
\vskip0.5pc
{\bf Proof} [Case 1] If the hereditary family ${\cal F}_h \cap [M]^{< \omega}$ is not pointwise closed,
then there exists $L \in [M]$ such that 
$[L]^{< \omega} \subseteq {\cal F}_h \subseteq {\cal F}$,  according to Proposition 3.6.
\par
[Case 2] Let ${\cal F}_h \cap [M]^{< \omega}$ be pointwise closed. Then
$\zeta^{\cal F}_M$ is a countable ordinal, according to Theorem 3.7. 

2(i) Let $\xi +1 < \zeta^{\cal F}_M$. Then $\xi +1 < \xi^{{\cal F}_h}_M$. 
According to Theorem 3.7 (subcase 2(i)) there exists $L \in [M]$ such that 
${\cal A}_{\xi} \cap [L]^{< \omega} \subseteq {\cal F}_h \subseteq {\cal F} \ . $
\par
2(ii) Let $\zeta^{\cal F}_M < \xi$ and $I \in [M]$. Then,
according to Theorem 3.7 (subcase 2(ii)) 
there exists $M_1 \in [I]$ such that
$$ {\cal A}_{\zeta^{\cal F}_M } \cap [M_1]^{< \omega} \subseteq [\Bbb N]^{< \omega} \setminus {\cal F}_h  
\ . \leqno (\star) $$

Using the Ramsey partition theorem for the countable ordinal $\xi$  (Theorem 1.5),  there exists
 an infinite subset 
$L$ of $M_1$ such that

either ${\cal A}_{\xi} \cap [L]^{< \omega} \subseteq {\cal F}$,
or $ {\cal A}_{\xi} \cap [L]^{< \omega} \subseteq [\Bbb N]^{< \omega} 
   \setminus {\cal F}$.
\par\noindent
We claim that the first alternative does not hold. 
Indeed, let ${\cal A}_{\xi}\cap [L]^{< \omega} \subseteq {\cal F}$.
Then $({\cal A}_{\xi} \cap [L]^{< \omega})^{\star} \subseteq {\cal F}^{\star} = {\cal F}$.
Using the canonical representation of  every infinite subset of $\Bbb N$ with respect to ${\cal A}_{\xi}$ (Proposition 2.4), it is easy to check that
\par\noindent
$({\cal A}_{\xi})^{\star} \cap [L]^{< \omega} = ( {\cal A}_{\xi} \cap [L]^{< \omega})^{\star} \ .$
Hence, $({\cal A}_{\xi})^{\star} \cap [L]^{< \omega}\subseteq {\cal F}$.

Since $\xi > \zeta^{{\cal F}}_M$ and  according to Corollary 3.2, there exists $L_1 \in [L]$ such that
$ ({\cal A}_{\zeta^{\cal F}_M })_{\star} \cap [L_1]^{< \omega} \subseteq 
({\cal A}_{\xi})^{\star} \cap [L]^{< \omega} \subseteq {\cal F}$
 and consequently 
  $({\cal A}_{\zeta^{\cal F}_M })_{\star} \cap [L_1]^{< \omega} \subseteq {\cal F}_h$.
This is a contradiction to $(\star)$; hence,
$ {\cal A}_{\xi} \cap [L]^{< \omega} \subseteq [\Bbb N]^{< \omega} 
   \setminus {\cal F} $ and  equivalently
 ${\cal F} \cap [L]^{< \omega} \subseteq ({\cal A}_{\xi})^{\star}
   \setminus {\cal A}_{\xi} \ , $ according to Proposition 3.1

2(iii) In the cases $\zeta^{\cal F}_M = \xi +1$ or $\zeta^{\cal F}_M = \xi$ we use Theorem 1.5.
\vskip1pc
{\bf Corollary 3.11 (=Theorem C,  Stronger form of Nash-Williams theorem in Gowers reformulation)} 
Let $\cal F$ be a tree of finite subsets of $\Bbb N$. Then there exists an infinite subset $L$ of 
$\Bbb N$, such that 

either (i) $[L]^{< \omega} \subseteq {\cal F}$;
\par
or  \ \     (ii)  there is a countable ordinal $\xi_0$, such that for every infinite subsets 
$I$ of $L$, there exists an initial segment $s$ of $I$ which belongs to $[\Bbb N]^{< \omega} \setminus {\cal F}$,
  and which is that unique initial segment of $I$ that belongs to ${\cal A}_{\xi_0}$.
\vskip1pc
{\bf Proof}  We apply  Theorem 3.10 (=Theorem B$^{\prime}$ in $\cal F$). 

If [Case 1] of Theorem 3.10 holds, then there exists 
$L \in [\Bbb N]$ such that $[L]^{< \omega} \subseteq {\cal F}$.

If [Case 2] of Theorem 3.10 holds, then there is a countable ordinal $\xi_0 = \zeta^{\cal F}_{\Bbb N}+1$
and $L \in [\Bbb N]$ such that 
${\cal A}_{\xi_0} \cap [L]^{< \omega} \subseteq [\Bbb N]^{< \omega} \setminus {\cal F}$.
According to Proposition 2.4 every infinite subset $I$ of $L$, has unique canonical representation with respect to 
${\cal A}_{\xi_0}$, hence for every $I \in [L]$ there exists a unique initial segment 
$s_{\xi_0,I}$ of I that belongs to ${\cal A}_{\xi_0}$ and consequently to $[\Bbb N]^{< \omega} \setminus {\cal F}$.
\vskip1pc
{\bf Remark 3.12}  Theorem C is indeed a stronger form of the classical Nash-Williams partition theorem, because it implies the Gowers reformulation of the Nash-Williams partition theorem (as given in the introduction of this paper).
To see that indeed Theorem C implies Gowers reformulation, let $\cal F$ be any family of finite subsets of $\Bbb N$. We set 
\par 
${\cal F}_t = \{ s \in {\cal F} $ : every non empty initial segment of $s$ belongs to 
${\cal F} \} \cup \{ \emptyset \}$.

The family ${\cal F}_t$ is a tree contained in $\cal F$. We apply Theorem C on ${\cal F}_t$. It follows that there exists an infinite subset $L$ of $\Bbb N$ such that.

either (i)$[L]^{< \omega} \subseteq {\cal F}_t$ (and consequently $[L]^{< \omega} \subseteq {\cal F}$)

or \ \ \ (ii) for every infinite subset $I$ of $\Bbb N$ there exists an initial segment $u$ of $I$ which
belongs to  
$[\Bbb N]^{< \omega} \setminus {\cal F}_t = ({\cal F} \setminus {\cal F}_t) \cup
 ([\Bbb N]^{< \omega} \setminus {\cal F}$. Thus 
 either $u \in [\Bbb N]^{< \omega} \setminus {\cal F}$ (in which case we set $s=u$), or 
 $u \in {\cal F} \setminus {\cal F}_t$ (in which case, by the definition of ${\cal F}_t$, there is a non empty initial segment $s$ of $u$ so that $s \in [\Bbb N]^{< \omega} \setminus {\cal F}$). Hence, in any of the two cases in (ii), for every infinite subset $I$ of $\Bbb N$, there is an initial segment $s$ of $I$ which belongs to
$[\Bbb N]^{< \omega} \setminus {\cal F}$, proving the Gowers reformulation of Nash-Williams theorem.
\vskip1pc
{\bf Remark 3.13} Gowers notices in [G], that if the first alternative (i) of the reformation of Nash-Williams's 
theorem does not hold, then
$[\Bbb N]^{< \omega} \setminus {\cal F}$ is large in an obvious sense and Nash-Williams's theorem asserts that if 
$[\Bbb N]^{< \omega} \setminus \cal F $ is a large subset of $[\Bbb N]^{< \omega}$, then there is an infinite subset $L$ of $\Bbb N$ for which $[\Bbb N]^{<\omega} \setminus \cal F $ has a stronger largeness property 
(alternative (ii)).
Theorem 3.10 (=Theorem B$^{\prime}$) is stronger than the Nash-Williams's theorem in the part that in the second alternative (ii)  the initial segments are located (uniformly for all infinite subsets) in the family 
${\cal A}_{\xi}$ and consequently $[\Bbb N]^{< \omega} \setminus {\cal F}$
has a much stronger largeness property than the given by  Nash-Williams's theorem.
\vskip0.5pc
Finally we state  our strengthening of the Nash-Williamsn [N-W] partition theorem in its original formulation
 the one concerning of pointwise closed families of infinite subsets of $\Bbb N$.

Firstly, we will give the necessary definitions.

\vskip1pc
{\bf Definition 3.14}
Let  $M$ be an infinite subset of $\Bbb N$, $s$ a finite subset of $\Bbb N$ and $\xi$ a countable ordinal. We set

(i) $[s,M]=\{s \cup L : L \in [M]$ and $s < L\}$,  \ \ \ $[\emptyset, M] = [M]$.

(ii) $s_{\xi,M}$ is the unique initial segment of $M$ which is an element of ${\cal A}_{\xi}$
(according to Proposition 2.4); note that  $ s_{0,M}= \emptyset$.
\vskip1pc
{\bf Theorem 3.15 (Stronger form of Nash-Williams's theorem)}
Let ${\cal U}$ be a pointwise closed family of infinite subsets of $\Bbb N$  and $M$ an infinite subset 
of $\Bbb N$. Then 

either (i) there exists $L \in [M]$ such that $[L] \subseteq {{\cal U}}$;

or \ \ \ (ii) there exists a countable ordinal $\zeta^{\cal U}_M$ such that for every 
countable ordinal   $\xi$ with 
$\xi > \zeta^{\cal U}_M$ and every $M_1 \in [M]$ there exists
$L \in [M_1]$ 
such that for every infinite subset $I$ of $L$ the unique initial segment 
$s_{\xi,I}$ of $I$ that belongs to ${\cal A}_{\xi}$ satisfies the relation
 $[s_{\xi, I}, \Bbb N] \subseteq [\Bbb N] \setminus {\cal U}$.
\vskip 1pc
{\bf Proof}
Let ${\cal F} = \{ s \in [\Bbb N]^{< \omega}$: $[s, \Bbb N] \cap {\cal U} \ne \emptyset\}$.
Of course ${\cal F}$ is a tree. We use Theorem 3.10.

If [Case 1] of Theorem 3.10 holds, then
 there exists $L \in [M]$ such that
$[L]^{< \omega} \subseteq {\cal F}$.
Then, $[s, \Bbb N] \cap {\cal U} \ne \emptyset$ for every $s \in [L]^{< \omega}$.
This gives that $[L] \subseteq {\cal U}$, since ${\cal U}$ is a pointwise closed family.

If [Case 2] of Theorem 3.10 holds, then  setting 
$\zeta^{\cal U}_M = \zeta^{\cal F}_M$ we have $\zeta^{\cal U}_M < \omega_1$ and for every
$\xi > \zeta^{\cal U}_M$ and every $M_1 \in [M]$ there exists $L \in [M_1]$ such that 
${\cal A}_{\xi} \cap [L]^{< \omega} \subseteq [\Bbb N]^{< \omega} \setminus {\cal F}$.
For every $I \in [L]$ let $s_{\xi,I}$ be the unique initial segment of $I$ which is an element of ${\cal A}_{\xi}$
(Proposition 2.4). Then $ s_{\xi,I} \in [\Bbb N]^{< \omega}\setminus {\cal F}$ for every
$I \in [\L]$. Hence, $[s_{\xi,I}, \Bbb N] \subseteq [\Bbb N]\setminus {\cal U}$ for every $I \in [L]$.
\vskip0.5pc
Immediate consequence of Theorem 3.15 is the classical Nash-Williams partition theorem:
\vskip1pc
{\bf{Corollary  3.15 (Nash-Williams [N-W])}}
Let $\cal U$ be a pointwise  closed family of infinite subsets of $\Bbb N$ and $M$ an infinite subset of $\Bbb N$.
Then 

either (i) there exists $L \in [M]$ such that $[L] \subseteq \cal U$; 

or \ \ \ (ii) there exists $L \in [M]$ such that $[L] \subseteq [\Bbb N] \setminus \cal U$,
equivalently, such that for every infinite subset $I$ of $L$ there exists an initial segment $s$ of $I$ such that $[s, \Bbb N] \subseteq \setminus {\cal U}$.
\vskip2pc
{\lfont 4. The derivation of Ellentuck's theorem}
\vskip2pc
We finally show that our Theorem 3.10(= Theorem B$^{\prime}$) implies, using the simple argument contained in 
Theorem 4.6, Ellentuck's theorem (and hence, Galvin-Prikry's and Silver's).

We recall the definition of the completely Ramsey families, given initially in [G,P] and [S].
\vskip1pc
{\bf Definition 4.1} A family $\cal U$ of infinite subsets of $\Bbb N$ is called {\bf completely Ramsey} if for every 
$\alpha \in [\Bbb N]^{< \omega}$ and $M \in [\Bbb N]$ there exists 
$L \in [M]$ such that

either (i) $[\alpha, L] \subseteq {\cal U}$;

or \ \ \ (ii) $[\alpha, L] \subseteq [\Bbb N] \setminus {\cal U}$.
\vskip1pc
{\bf Theorem 4.2}
Let ${\cal U}$ be a pointwise close  family of infinite subsets of $\Bbb N$, $\alpha$ a finite subset of $\Bbb N$ with 
cardinality $m$  and $M$ an infinite subset of $\Bbb N$. Then 

either (i) there exists $L \in [M]$ such that $[\alpha, L] \subseteq {\cal U}$;

or \ \ \ (ii) there exists a countable ordinal $\zeta^{\cal U}_{\alpha,M}$ such that for every 
countable ordinal $\xi$ with  
$\xi > \zeta^{\cal U}_M$ and every $M_1 \in [M]$  there exists
$L \in [M_1]$ such that \par\noindent
$[\alpha \cup s_{\xi,I}, \Bbb N] \subseteq [\Bbb N] \setminus {\cal F}$ \ for every infinite subset
$I$ of $L$.
\vskip 1pc
{\bf Proof}
Let ${\cal F} = \{ s \in [\Bbb N]^{< \omega}$: $[\alpha \cup s, \Bbb N] \cap {\cal U} \ne \emptyset\}$.
Of course ${\cal F}$ is a tree. We use Theorem 2.10.
\vskip1pc
{\bf Corollary 4.3 (Galvin-Prikry [G-P])} Every pointwise closed (resp. pointwise open) family of infinite subsets of $\Bbb N$ is completely Ramsey.
\vskip1pc
{\bf Definition 4.4}
Ellentuck's topology on $[\Bbb N]$ is the topology which has base the family of all sets
$[\alpha, M]$, where $\alpha \in [\Bbb N]^{< \omega}$ and $M \in [\Bbb N]$.
Of course Ellentuck's topology is weaker than the topology of pointwise convergence.
\par\noindent
We denote by $\hat{\cal U}$ and $\cal{U}^{\diamond}$ the closure and the interior respectively of a subset 
$\cal{U}$ of $\Bbb N$ in the Ellentuck's topology. Then, it is easy to see that 

(i) $\hat {\cal U} = \{ L \in [\Bbb N] : [s_{n,L}, L] \cap {\cal U} \ne \emptyset$ for every $n \in \Bbb N$;
and 

(ii) ${\cal U}^{\diamond}= \{ L \in [\Bbb N] :$ there exists $n \in \Bbb N$ such that
$[s_{n,L}, L] \subseteq {\cal U} \}$ .
\vskip1pc
{\bf Lemma 4.5}
Let ${\cal L} \subseteq \{[s,I]: \ \ s \in [\Bbb N]^{<\omega} \ {\rm and} \ I \in [\Bbb N]\}$ with the following two properties:

(i) For every $(\alpha, M) \in [\Bbb N]^{< \omega} \times [\Bbb N]$ 
there exists $I \in [M]$ such that $[\alpha, I] \in {\cal L}$; and

(ii) if $[s, I] \in {\cal L}$, then $[s, L] \in {\cal L}$ for every $L \in [I]$.
\par\noindent
Then for every $(\alpha, M) \in [\Bbb N]^{< \omega} \times [\Bbb N]$ there exists 
$L \in [\alpha, M]$ such that $[\alpha \cup \beta, L] \in {\cal L}$ for every 
$\beta \in [L]^{< \omega}$.
\vskip1pc
{\bf Proof}
Let $(\alpha, M) \in [\Bbb N]^{< \omega} \times [\Bbb N]$ and let $m$ be the cardinality of $\alpha$. 
We can assume that 
$\alpha \prec M$. Set $L_0 = M$. According to property (i) of ${\cal L}$ there exists $L_1 \in [M]$ such that 
$[\alpha, L_1] \in {\cal L}$. Set $L_1 = I$. Let $L_1 \subseteq \ldots \subseteq L_n$ have been constructed and let
$\{s_1, \ldots , s_r \} = \{ s \in [L_n]^{< \omega}$ : $s \subseteq s_{m+n,L_n}$ and  
$s \not\subseteq s_{m+n-1, L_n} \}$. According to property (i) of $\cal L$ there exists 
$I_{n+1}^1 \in [L_n]$ such that $[s_1, I^1_{n+1}] \in {\cal L}$.
Setting 
\par\noindent
$L^1_{n+1} = I_{n+1}^1 \cup s_{m+n, L_n}$ we have that $L^1_{n+1} \in [s_{m+n, L_n}, L_n]$ and that
$[s_1, L^1_{n+1}] = [s_1, I^1_{n+1}] \in {\cal L}$. Analogously, we can choose 
$L^2_{n+1} \in [s_{m+n, L_n}, L^1_{n+1}]$
such that $[s_2, L^2_{n+1}] \in {\cal L}$ and so on. Set
$L_{n+1}=L^r_{n+1}$.

Since $L_{n+1} \in [s_{m+n,L_n}, L_n]$ for every $n \in \Bbb N$ there exists $L \in [M]$ such that
$s_{m+n, L} = s_{m+n, L_n}$ for every $n \in \Bbb N$.  Hence, $L$ 
has the desired property, according to property (ii) of $\cal L$. 
\vskip1pc
{\bf Theorem 4.6}
Let $\cal U$ be a family of infinite subsets of $\Bbb N$, $M$ an
 infinite subset of $\Bbb N$ and $\alpha$ a finite subset  of $\Bbb N$. Then 
there exists $L \in [M]$ such that

either (i)
$[\alpha, L] \subseteq \hat {\cal U}$ ;

or \ \ \ (ii) 
$ [\alpha, L] \subseteq  [\Bbb N] \setminus {\cal U}$.
\vskip1pc
{\bf Proof}
Let $(\alpha, M) \in [\Bbb N]^{< \omega} \times [\Bbb N]$.  Set
\par\noindent
${\cal L}^{\cal U} = \{ [s, I]:$ 
either $[s, I] \cap {\cal U} = \emptyset$
\ \ or $[s, I_1] \cap {\cal U} \ne \emptyset$ for every $I_1 \in [I] \}$.
\par\noindent
It is easy, to check that ${\cal L}^{\cal U}$ satisfies the assumptions (i) and (ii) of Lemma 3.9, 
hence there exists $I \in [\alpha, M]$ such that 
 $ [\alpha \cup \beta, I] \in {\cal L}^{\cal U}$ for every $\beta \in [I]^{< \omega}$.
We assume that $[\alpha, I_1] \cap {\cal U} \ne \emptyset$ for every $I_1 \in [I]$.

Set 
${\cal F} = \{ \beta \in [I]^{< \omega} :$ $\alpha < \beta $ and 
$[\alpha \cup \beta, I_1] \cap {\cal U} \ne \emptyset$ for every
$I_1 \in [I]\}$.
The family ${\cal F} $ is a tree.  We use  Theorem 3.10. 

If [Case 1] of Theorem 3.10 holds, then there exists 
$L \in [I]$ such that $[L]^{< \omega} \subseteq {\cal F}$. Then $[\alpha, L] \subseteq \hat {\cal U}$.  

[Case 2] of Theorem 3.10 does not occur. Let  
${\cal A}_{\xi} \cap [L]^{<\omega} \subseteq [\Bbb N]^{< \omega} \setminus {\cal F}$ for some $\xi < \omega_1$.
Then $[\alpha, L] \cap {\cal U} = \emptyset$. Indeed, let $L_2 \in [\alpha, L] \cap \cal U$, 
 $L_2 = \alpha \cup L_1, L_1 \in [L]$ and $\alpha < L_1$. Then $L_2 \in [\alpha \cup s_{\xi,L_1}, I ] \cap {\cal U}$
and consequently $s_{\xi,L_1} \in {\cal A}_{\xi} \cap [L]^{< \omega} \cap {\cal F}$.
A contradiction; hence $[\alpha, L] \cap {\cal U} = \emptyset$.
This is a contradiction to our assumption that $[\alpha, I_1] \cap {\cal U} \ne \emptyset$
for every $I_1 \in [I]$.

Hence, either there exists $L \in [M]$ such that $ [\alpha, L] \subseteq  [\Bbb N] \setminus {\cal U}$ 
or there exists $L \in [M]$ such that $[\alpha, L] \subseteq \hat {\cal U}$. 
\vskip1pc
{\bf Corollary 4.7}
Every family of infinite subsets of $\Bbb N$ which is closed (resp. is open) in the Ellentuck's topology,  is completely Ramsey. 
\vskip1pc
{\bf Corollary  4.8}
Let $\cal U$ be a family of infinite subset of $\Bbb N$ which is a meager set in the Ellentuck's topology, $M \in [\Bbb N]$ and 
$\alpha \in [\Bbb N]^{< \omega}$. Then there exists $L \in [M]$ such that
$ [\alpha, L] \subseteq [\Bbb N] \setminus {\cal U}$.
\vskip1pc
{\bf Proof}
Let ${\cal U}= {\displaystyle \bigcup^{\infty}_{n=0}} {\cal U}_n$ where $(\hat{\cal U}_n)^{\diamond} = \emptyset$
for every $n = 0,1, \ldots$. According to Theorem 3.10, there exists $L \in [M]$ such that

either (i)  $[\alpha, L] \subseteq \hat {\cal U}$; 

or \ \ \ (ii) $[\alpha, L]   \subseteq [\Bbb N] \setminus {\cal U}$.
\par\noindent
We will prove that the first alternative is impossible. Let $L \in [M]$. Set 
\par\noindent
${\cal L} = \{[s,I]: s \in [\Bbb N]^{< \omega}, I \in [\Bbb N]$ and
$[s,I] \cap {\cal U}_{\kappa} = \emptyset$  for every  $\kappa \in \Bbb N$ with 
$\kappa \leq \vert s \vert \}$.
\par\noindent
The family $\cal L$ satisfies the assumptions (i) and (ii) of Lemma 3.9, according to Theorem 3.10.
 Hence there exists 
$L_1 \in [\alpha, L]$ such that
$[ \alpha \cup \beta, L_1] \cap {\cal U}_{\kappa} = \emptyset$ for every $\beta \in [L_1]^{< \omega}$ and 
 $\kappa \in \Bbb N$ with $\kappa \leq \vert \alpha \cup \beta \vert$.
Then $L_1 \not\in \hat {\cal U}$, since $[\alpha, L_1] \cap {\cal U} = \emptyset$.
Indeed, if $L_2 \in [\alpha, L_1]\cap {\cal U}$, then $L_2\in {\cal U}_{\kappa}$ for some 
$\kappa =0,1, \ldots $ and choosing 
an initial segment $\beta$ of $L_2$ such that $\vert \alpha \cup \beta \vert \geq \kappa$
we have $L_2 \in [\alpha \cup \beta, L_2] \cap {\cal U}_{\kappa} = \emptyset$. A contradiction; hence 
$I_1 \not\in \hat{\cal U}$ and consequently $[\alpha, L] \not\subseteq \hat{\cal U}$.
\vskip1pc
{\bf Corollary 4.9} {\bf (Ellentuck [E])}
A family $\cal U$ of infinite subsets of $\Bbb N$ is completely Ramsey if and only if $\cal U$
has the Baire property in Ellentuck's topology.
\vskip1pc 
{\bf Proof}
Let ${\cal U}$ has the Baire property in Ellentuck's topology. Then, setting 
${\cal C}^c = [\Bbb N] \setminus{\cal C}$ for every ${\cal C} \subseteq [\Bbb N]$, we have 
${\cal U} = {\cal B} \triangle {\cal C}= ({\cal B} \cap {\cal C}^c)\cup ({\cal C} \cap {\cal B}^c)$
 where $\cal B$ is a closed  set and $\cal C$  a meager set in Ellentuck's topology. According  to Corollary 3.12, there exists $L_1 \in [\Bbb N]$ such that 
$[\alpha, L_1] \subseteq {\cal C}^c$. According to Theorem 3.10, there exists $L \in [L_1]$ such that

either (i) $[\alpha, L] \subseteq {\cal B} \cap {\cal C}^c \subseteq {\cal U}$;

or \ \ \ (ii) $[\alpha, L] \subseteq {\cal B}^c \cap {\cal C}^c \subseteq [\Bbb N] \setminus {\cal U}$
\par\noindent
Hence $\cal U$ is completely Ramsey.

On the other hand, if ${\cal U}$ is completely Ramsey, then ${\cal U}$ has the Baire property in Ellentuck's topology, since  ${\cal U} =  {\cal U}^{\diamond} \cup ({\cal U} \setminus {\cal U}^{\diamond})$ and
${\cal U} \setminus {\cal U}^{\diamond}$ is a meager set in Ellentuck's topology.
\vskip1pc
{\bf Remark 4.10 (i) (Galvin-Prikry [G-P])}
Every  family of finite subsets of $\Bbb N$ which is a Borel set  in the topology of pointwise convergence is completely Ramsey,
since every Borel set has the Baire property.

{\bf (ii)  (Silver [S])}
Every family of finite subsets of $\Bbb N$ which is an analytic set in the topology of pointwise convergence 
 is completely Ramsey, since every analytic set has the Baire property.
\vskip2pc\noindent
{\centerline{\lfont Bibliography}}
\vskip2pc\noindent
{\parindent=1.5cm
\litem{[A-A]}
D. Alspach and S. Argyros, Complexity of weakly null
sequences, Dissertations Math. 321 (1992), 1--44.
\litem{[A-M-T]}
S. Argyros, S. Mercourakis and A. Tsarpalias, Convex
unconditionality and summability of weakly null sequences, Israel Journal of
Math. 107 (1998), 157-193.
\litem{[A-O]}
D. Alspach and E. Odell, Averaging weakly null
sequences, Lecture Notes in Math. 1332, Springer, Berlin, 1988.
\litem{[B]}
I. Bendixson, Quelques theor\`{e}mes de la th\'{e}orie des
ensembles de points, Acta Math. 2 (1883), 415--429.
\litem{[C1]} 
G. Cantor, Grundlagen einer allgemeine
Mannigfaltigkeitslehre, Math. Annalen, 21 (1883), 575.
\litem{[C2] } 
G. Cantor, Beitr\"age zur Begr\"undung der transfiniten
Mengenlehre II, Math. Ann. 49, (1897) 207 - 246.
\litem{[E]} 
E.E. Ellentuck, A new proof that analytic sets are Ramsey,
J. Symbolic Logic 39 (1974), 163--165.
\litem{[F1]} 
V. Farmaki, On Baire-1/4 functions and spreading models, Mathemateka, 41(1994), 251-265.
\litem{[F2]} 
V. Farmaki, Classifications of Baire--1 functions and
$c_0$--spreading models, Trans. Amer. Math. Soc. 345 (2), (1994), 819--831.
\litem{[F3]}
V. Farmaki, On Baire-1/4 functions, Trans.Amer. Math. Soc. 348(10), (1996), 4023-4041.
\litem{[F4]} 
V. Farmaki, Ramsey dichotomies with ordinal index, arXiv: math. LO/9804063 v1, 1998, electronic prepublication.
\litem{[F5]} 
V. Farmaki, The uniform convergence ordinal index and the $l^1$-behavior of a
sequence of functions, Positivity . 
\litem{[F6]} 
V. Farmaki, Ordinal indices and Ramsey dichotomies measuring $c_0$-content and
semi bounded completeness, Fundamenta Mathematicae 172(2002) 153-179
\litem{[G-P]} 
F. Galvin and K. Prikry, Borel sets and Ramsey's
theorem, J. Symbolic Logic 38 (1973), 193--198.
\litem{[G]} 
W.T. Gowers, An infinite Ramsey theorem and some Banach-space dichotomies, Annals of Mathematics, 156 (2002), 797-833.
\litem{[J]} 
R. Judd, A dichotomy on Schreier sets, Studia Math. 132 (1999), 245-256.
\litem{[K]} 
K. Kuratowski, Topology, Valume I, Academic press, (1966).
\litem{[L] } 
A. Levy, Basic set Theory, Springer-Verlag, (1979).
\litem{[M-N]} 
S. Mercourakis and S. Negrepontis, Banach spaces and
Topology II, Recent Progress in
General Topology, M. Husek and J. Vaan Mill (editors) Elsevier Sciences
Publishers, (1992).
\litem{[N-W]} 
C.St.J.A. Nash--Williams, On well quasi-ordering
transfinite sequences, Proc. Camb. Phill. Soc. 61 (1965), 33--39.
\litem{[O]} 
E. Odell, Applications of Ramsey theorems to Banach space
theory, Notes in Banach spaces (H.E. Lacey, ed.), Univ. of Texas Press,
(1980), 379--404.
\litem{[O-T-W]} 
E. Odell, N. Tomczak--Jaegermann and R. Wagner,
Proxinity to $\ell_1$ and distotion in asymptotic $\ell_1$ spaces,
(preprint). 
\litem{[P-R]} 
P. Pudl\'{a}k and V. R\"{o}dl, Partition theorems for
systems of finite subsets of integers, Discrete Math. 39 (1982), 67--73.
\litem{[R]} 
F.P. Ramsey, On a problem of formal logic, Proc. London
Math. Soc. 30(2), (1929), 264--286.
\litem{[RO1]} 
H. Rosenthal, Weakly independent sequences and the
Banach--Saks property, Bulletin London Math. Soc., 8 (1976), 22--24.
\litem{[RO2]} 
H. Rosenthal, A characterization of Banach spaces containing $c_0$, Journal of the Amer. Math. Soc.
7(3), (1994), 707-748.
\litem{[S]} 
J. Schreier, Ein Gegenbeispiel zur Theorie der schwachen Konvergenz,
Studia Math. 2(1930), 58-62.
\litem{[Si]} 
J. Silver, Every analytic set is Ramsey, J. Symbolic
Logic 35 (1970), 60--64. 
\litem{[TJ]} 
N. Tomczak-Jaegermann, Banach spaces of type p have  arbitrarily distortable subspaces, Geom. Funct.Anal.
6(1996),1074-1082.
\litem{[T]} 
B.S. Tsirelson, Non every Banach space contains $l^{p}$ or $c_0$, Funct. Annal. Appl. 8 (1974), 138-141.
\par}

\vskip2pc
Address of the author:
\par
Department of Mathematics,
\par
University of Athens,
\par
Panepistimiopolis,
157 84 Athens, Greece
\par
{e-mail: {vfarmaki@math.uoa.gr}
\par
\par
March 5, 2004

\end